\documentclass[12pt,reqno]{amsart}

\usepackage{amsmath,amssymb,amsfonts,amsthm}
\usepackage{palatino}
\usepackage[euler-digits]{eulervm}
\usepackage{microtype}
\usepackage{tikz,tikz-cd}
\usepackage{graphicx}
\usepackage{fullpage}
\usepackage{ytableau}
\usepackage{fancyhdr}
\usepackage{placeins} 

\usepackage{academicons}
\usepackage{xcolor}

\usepackage[margin=1in]{geometry}

\usepackage[
pdfauthor={},
pdfkeywords={},
pdftitle={},
pdfcreator={},
pdfproducer={},
linktocpage,colorlinks,bookmarksnumbered,linkcolor=blue,
citecolor=red,urlcolor=red]{hyperref}

\DeclareMathOperator{\Li}{Li}

\newcommand{\QED}{\hspace{\stretch{1}} $\blacksquare$}

\renewcommand{\(}{\left\(}
\renewcommand{\)}{\right\)}

\theoremstyle{plain}
\newtheorem{thm}{Theorem}

\newtheorem{cor}[thm]{Corollary}
\newtheorem{prop}[thm]{Proposition}

\newtheorem{claim}[thm]{Claim}

\theoremstyle{definition}

\theoremstyle{remark}

\numberwithin{equation}{section}
\numberwithin{thm}{section}
\allowdisplaybreaks

\definecolor{lime}{HTML}{A6CE39}

\newcommand{\orcid}[1]{\href{https://orcid.org/#1}{\textcolor[HTML]{A6CE39}{\aiOrcid}}}

\begin{document}
	\title[Evaluation of Multivariate Integrals]{Evaluation of Multivariate Integrals based on a Duality Identity for the Stieltjes Transform}
	\author{Sarth Chavan}
	\author{Christophe Vignat}
	\date{\today}  
	\thanks{The first author is supported by The 2022 Spirit of Ramanujan Fellowship.}
	\address{Euler Circle, Palo Alto, California 94306, USA}
	\email{sarth5002@outlook.com}
	\address{Department of Mathematics, Tulane University, New Orleans, USA; \newline .\quad Universit\'e Paris-Saclay, CentraleSupelec, Orsay, France}
	\email{cvignat@tulane.edu}
	\maketitle
	\begin{abstract}  
		A detailed study of a double integral representation of the Catalan's
		constant allows us to identify a duality identity for the Stieltjes
		transform on which it is based. This duality identity is then extended
		to an arbitrary dimensional integral and several special cases are
		deduced. On the way, we also highlight a relationship with some multivariate
		generalizations of the Riemann zeta function.
	\end{abstract}
	
	\maketitle
	
	\section{Introduction}
	
	This study started as the first author, looking at integral representations
	for the Catalan's constant 
	\[
	G=\sum_{n=0}^{\infty}\dfrac{\left(-1\right)^{n}}{\left(2n+1\right)^{2}},
	\]
	computed the value of the integral 
	\begin{equation}
		\int_{0}^{\pi/2}\int_{0}^{\pi/2}\dfrac{\log\left(\cos\left(\dfrac{x}{2}\right)\right)-\log\left(\cos\left(\dfrac{z}{2}\right)\right)}{\cos x-\cos z}\,\mathrm{d}x\,\mathrm{d}z=\pi G-\frac{7}{4}\,\zeta(3),
		\label{double_integral}
	\end{equation}
	to deduce the following representation.
	\begin{prop}
		\label{Prop1} The Catalan constant is equal to 
		\begin{equation}
			G=\dfrac{1}{\pi}\left\{ \dfrac{7}{4}\,\zeta(3)+\int_{0}^{\pi/2}\int_{0}^{\pi/2}\dfrac{\log\left(\cos\left(\dfrac{x}{2}\right)\right)-\log\left(\cos\left(\dfrac{z}{2}\right)\right)}{\cos x-\cos z}\,\mathrm{d}x\,\mathrm{d}z\right\}. \label{G}
		\end{equation}
	\end{prop}
	
	Notice that this double integral \eqref{double_integral} appears in
	Problem 1.53 in \cite{Valean} and its computation in \cite[Solution 3.53]{Valean}
	produces identity \eqref{G}. 
	
	Our first attempts at generalizing this result produced variations such
	as 
	\begin{prop}
		\label{cosh_int} The following identities hold
		
		\[
		\int_{0}^{\pi/2}\int_{0}^{\pi/2}\dfrac{\log\left(\cos x\right)-\log\left(\cos z \right)}{\cos x-\cos z}\,\mathrm{d}x\,\mathrm{d}z =2\pi G,
		\]
		\[
		\int_{0}^{\infty}\int_{0}^{\infty}\frac{\log\left(\cosh\left(\dfrac{x}{2}\right)\right)-\log\left(\cosh\left(\dfrac{z}{2}\right)\right)}{\cosh x-\cosh z}\,\mathrm{d}x\,\mathrm{d}z=\dfrac{7}{4}\,\zeta(3),
		\]
		\[
		\int_{0}^{\pi/2}\int_{0}^{\pi/2}\frac{\log\left(1+\tan\vartheta_{1}\right)-\log\left(1+\tan\vartheta_{1}\right)}{\tan\vartheta_{1}-\tan\vartheta_{2}}\mathrm{d}\vartheta_{1}\mathrm{d}\vartheta_{2}=\frac{\pi}{16}\left(\pi^{2}+\pi-4\log2\right) -G.
		\]
	\end{prop}
	The proofs of the first two identities appear in Subsection \ref{subsec:Proof Cosh} and for the last one, in Subsection \ref{proof of proposition 8.3}. However, in order to obtain more general integrals and relate them to other interesting quantities such as multivariate sums, more elaborate
	techniques are required. The aim of this paper is to exhibit some of these techniques. 
	\subsection{Organization of the manuscript} In Section \ref{sec:A-Fourier-series}, a Fourier
	expansion is used to deduce the value of the tail of a divergent multiple
	zeta value (MZV) of depth 3. Fourier expansion is also used in Section \ref{sec:A-Generating-Function}
	to obtain a generalization of the integral \eqref{double_integral}
	to other functions than the logarithmic one. In Section \ref{sec:A-Triple-Integral},
	it is shown that the double integral \eqref{double_integral} has an equivalent nested triple
	integral and nested triple sum (MZV) expression. Section \ref{sec:Integration-over}
	produces values of similar integrals when the integration domain is
	changed from $\left[0,\pi/2\right]^{2}$ to $\left[0,\pi\right]^{2}.$
	Section \ref{sec:A-parameterized-version} exhibits a parameterized generalization of integral \eqref{double_integral} from which its value
	is recovered for a special value of the parameter. Section \ref{sec:A-simple-integral}
	uses the notion of Hadamard product to obtain a one-dimensional integral version
	of a more general version of integral \eqref{double_integral}. Finally,
	Section \ref{sec:multivariate-extensions} derives multivariate generalizations
	of the same integral based on a duality result for analogues of the Stieltjes
	transform. The last section, Section \ref{proof}, produces the proofs of all the results obtained in this manuscript. 
	
	Beyond this specific case, we hope that the tools used in this study will    be useful to the readers in their endeavour of evaluating more interesting integrals.
	\section{\label{sec:A-Fourier-series}A Fourier Series Computation}
	This section provides a computation of the double integral \eqref{double_integral}
	based on the Fourier series expansion of $\log\left(2\cos\frac{\vartheta}{2}\right)$, see for example  \cite[Identity 1.441.4]{Gradshteyn}
	\begin{equation}
		\log\left(2\cos\frac{\vartheta}{2}\right)=\sum_{k=1}^{\infty}\dfrac{(-1)^{k-1}}{k}\cos\left(k\vartheta\right),\quad\left(-\,\pi < \vartheta < \pi\right),
		\label{eq:Fourier logcos}
	\end{equation}
	which produces an unexpected link with the multiple zeta values (MZVs). 
	\begin{prop}
		\label{prop:Prop2}From the value of the double integral 
		\[
		\int_{0}^{\pi/2}\int_{0}^{\pi/2}\dfrac{\log\left(\cos\left(\dfrac{x}{2}\right)\right)-\log\left(\cos\left(\dfrac{z}{2}\right)\right)}{\cos x-\cos z}\,\mathrm{d}x\,\mathrm{d}z=\pi G-\dfrac{7}{4}\,\zeta(3),
		\]
		we deduce the value of the unusual and nicely symmetrical nested triple sum 
		\begin{equation}
			\sum_{k=1}^{\infty}\frac{1}{2k+1}\sum_{\ell\geqslant k+1}\frac{\left(-1\right)^{\ell+1}}{2\ell}\sum_{m=0}^{\ell-1}\frac{1}{2m+1}=\frac{7}{32}\,\zeta(3)-\dfrac{3}{16}\,\zeta(2).\label{eq:MZV1}
		\end{equation}
	\end{prop}
	
	The proof of this proposition is given in Subsection \ref{ProofProp2} and generates, as by-products, the two evaluations (see \eqref{prud_identity} and
	\eqref{eq:e2n+1}):
	\begin{equation}
		\int_{0}^{\pi/2}\int_{0}^{\pi/2}\dfrac{\cos\left(2kx\right)-\cos\left(2kz\right)}{\cos x-\cos z}\,\mathrm{d}x\,\mathrm{d}z=2\pi\sum_{\ell=1}^{k-1}\dfrac{\left(-1\right)^{\ell}}{2\ell+1},\label{eq:by-product1}
	\end{equation}
	and
	\begin{equation}
		\int_{0}^{\pi/2}\int_{0}^{\pi/2}\dfrac{\cos\left((2k+1)x\right)-\cos\left((2k+1)z\right)}{\cos x-\cos z}\mathrm{d}x\mathrm{d}z=\frac{\pi^{2}}{4}+4\sum_{k=1}^{n}\frac{\left(-1\right)^{k+1}}{k}\sum_{m=0}^{n-1}\frac{1}{2m+1}.\label{eq:by-product2}
	\end{equation}
	Notice that the multiple zeta value
	\begin{equation}
		\sum_{k=1}^{\infty}\frac{1}{2k+1}\sum_{\ell=1}^{k}\frac{\left(-1\right)^{\ell+1}}{2\ell}\sum_{m=0}^{\ell-1}\frac{1}{2m+1},\label{eq:MZV1DV}
	\end{equation}
	is divergent: its companion form \eqref{eq:MZV1}, however, is finite
	and can be considered as the tail of the divergent series \eqref{eq:MZV1DV}.
	To our best knowledge, tails of divergent multiple zeta values and its analogues have not been studied in
	the literature.
	
	To conclude this section, let us check that identity \eqref{eq:MZV1}
	allows in return to recover the value of the double integral \eqref{double_integral}:
	first  notice that, using a simple integral representation for
	$\sum_{m=0}^{\ell-1}\frac{1}{2m+1},$ the value of the MZV \eqref{eq:MZV1}
	is
	\begin{align*}
		& \sum_{k=1}^{\infty}\dfrac{1}{2k+1}\sum_{\ell\geqslant k+1}^{\infty}\dfrac{\left(-1\right)^{\ell+1}}{\ell}\sum_{m=0}^{\ell-1}\dfrac{1}{2m+1}\\&=\sum_{k=1}^{\infty}\dfrac{1}{2k+1}\sum_{\ell\geqslant k+1}^{\infty}\dfrac{\left(-1\right)^{\ell+1}}{\ell}\int_{0}^{1}\dfrac{1-\vartheta_{1}^{2\ell}}{1-\vartheta_{1}^{2}}\,\mathrm{d}\vartheta_{1}
		\\&=\sum_{\ell=1}^{\infty}\dfrac{\left(-1\right)^{\ell+1}}{\ell}\sum_{k=1}^{\ell-1}\dfrac{1}{2k+1}\int_{0}^{1}\dfrac{1-\vartheta_{1}^{2\ell}}{1-\vartheta_{1}^{2}}\,\mathrm{d}\vartheta_{1}
		\\&=\sum_{l=1}^{\infty}\dfrac{\left(-1\right)^{\ell+1}}{\ell}\int_{0}^{1}\dfrac{1-\vartheta_{2}^{2\ell}}{1-\vartheta_{2}^{2}}\int_{0}^{1}\dfrac{1-\vartheta_{1}^{2\ell}}{1-\vartheta_{1}^{2}}\,\mathrm{d}\vartheta_{1}\,\mathrm{d}\vartheta_{2}-\sum_{\ell=1}^{\infty}\dfrac{\left(-1\right)^{\ell+1}}{\ell}\int_{0}^{1}\dfrac{1-\vartheta_{1}^{2\ell}}{1-\vartheta_{1}^{2}}\,\mathrm{d}\vartheta_{1}.
	\end{align*}
	Interchanging the sum and integral signs produces
	\begin{align*}
		&\int_{0}^{1}\dfrac{1}{1-\vartheta_{2}^{2}}\int_{0}^{1}\dfrac{1}{1-\vartheta_{1}^{2}} \sum_{\ell=1}^{\infty}\dfrac{(-1)^{\ell+1}(1-\vartheta_{1}^{2\ell})(1-\vartheta_{2}^{2\ell})}{\ell} \mathrm{d}\vartheta_{1}\mathrm{d}\vartheta_{2} -\int_{0}^{1}\dfrac{1}{1-\vartheta_{1}^{2}} \sum_{\ell=1}^{\infty}\dfrac{(-1)^{\ell+1}(1-\vartheta_{1}^{2\ell})}{\ell}\mathrm{d}\vartheta_{1}
		\\&=\int_{0}^{1}\int_{0}^{1}\dfrac{\left(\log\left(1+\vartheta_{2}^{2}\vartheta_{1}^{2}\right)-\log\left(1+\vartheta_{1}^{2}\right)-\log\left(1+\vartheta_{2}^{2}\right)+\log2\right)}{\left(1-\vartheta_{2}^{2}\right)\left(1-\vartheta_{1}^{2}\right)}\mathrm{d}\vartheta_{1}\mathrm{d}\vartheta_{2}+\int_{0}^{1}\dfrac{\log\left(\frac{1+\vartheta_{1}^{2}}{2}\right)}{1-\vartheta_{1}^{2}}\mathrm{d}\vartheta_{1}
		\\&=\int_{0}^{1}\int_{0}^{1}\dfrac{\vartheta_{1}\log\left(\frac{1-\vartheta_{1}}{1+\vartheta_{1}}\right)}{\left(1+\vartheta_{1}^{2}\right)\left(1+\vartheta_{2}^{2}\vartheta_{1}^{2}\right)}\,\mathrm{d}\vartheta_{1}\mathrm{d}\vartheta_{2}+\dfrac{1}{2}\int_{0}^{1}\dfrac{\log\left(1+\vartheta_{1}^{2}\right)}{\vartheta_{1}}\,\mathrm{d}\vartheta_{1}-\int_{0}^{1}\dfrac{\log\left(1+\vartheta_{1}\right)}{\vartheta_{1}}\mathrm{d}\vartheta_{1}.
	\end{align*}
	Some routine manipulations produce 
	\begin{align*}
		&\int_{0}^{1}\int_{0}^{1}\dfrac{\vartheta_{1}\log\left(\frac{1-\vartheta_{1}}{1+\vartheta_{1}}\right)}{\left(1+\vartheta_{1}^{2}\right)\left(1+\vartheta_{2}^{2}\vartheta_{1}^{2}\right)}\,\mathrm{d}\vartheta_{1}\,\mathrm{d}\vartheta_{2}+\dfrac{1}{2}\int_{0}^{1}\dfrac{\log\left(1+\vartheta_{1}^{2}\right)}{\vartheta_{1}}\,\mathrm{d}\vartheta_{1}-\int_{0}^{1}\dfrac{\log\left(1+\vartheta_{1}\right)}{\vartheta_{1}}\,\mathrm{d}\vartheta_{1}
		\\ & =-\dfrac{1}{4}\left\{ \pi\int_{0}^{1}\dfrac{\log\vartheta_{1}}{1+\vartheta_{1}^{2}}\,\mathrm{d}\vartheta_{1}+3\int_{0}^{1}\dfrac{\log\left(1+\vartheta_{1}\right)}{\vartheta_{1}}\,\mathrm{d}\vartheta_{1}\right\} +\int_{0}^{\pi/4}\vartheta_{1}\log\left(\tan\vartheta_{1}\right)\mathrm{d}\vartheta_{1}
		\\ & =\dfrac{\pi G}{4}+\dfrac{7}{16}\,\zeta(3)-\dfrac{\pi G}{4}-\dfrac{3}{8}\,\zeta(2)=\dfrac{7}{16}\,\zeta(3)-\dfrac{3}{8}\,\zeta(2),
	\end{align*}
	as desired, where we have used the integral evaluations
	\[
	\int_{0}^{1}\dfrac{\log\vartheta_{1}}{1+\vartheta_{1}^{2}}\mathrm{d}\vartheta_{1}=-G,\,\,\int_{0}^{1}\dfrac{\log\left(1+\vartheta_{1}\right)}{\vartheta_{1}}\mathrm{d}\vartheta_{1}=\frac{1}{2}\,\zeta(2),\,\,
	\int_{0}^{\pi/4}\vartheta_{1}\log\tan\vartheta_{1}\mathrm{d}\vartheta_{1}=\frac{7}{16}\,\zeta(3)-\frac{\pi G}{4}.
	\]
	\section{\label{sec:A-Generating-Function}A Generating Function Approach}
	
	In this section, the previous result is used to obtain a generalization
	of the original integral \eqref{double_integral} of the form 
	\begin{equation}
		\int_{0}^{\pi/2}\int_{0}^{\pi/2}\dfrac{\mathfrak{F}\left(\vartheta_{1}\right)-\mathfrak{F}\left(\vartheta_{2}\right)}{\cos\vartheta_{1}-\cos\vartheta_{2}}\,\mathrm{d}\vartheta_{1}\,\mathrm{d}\vartheta_{2},\label{eq:general integral}
	\end{equation}
	in the case where a Fourier series for the function $\mathfrak{F}$ is known. We assume in the rest of this section that the conditions for convergence and exchange of sums and integrals are met.
	\begin{thm}
		\label{thm:generating_function} Assume that the function $\mathfrak{F}$ has
		the Fourier series expansion 
		\begin{equation}
			\mathfrak{F}(\vartheta)=\sum_{k=0}^{\infty}\alpha_{k}\cos\left(k\vartheta\right),\label{eq:Fourier expansion}
		\end{equation}
		then the following identity holds
		\[
		\int_{0}^{\pi/2}\int_{0}^{\pi/2}\dfrac{\mathfrak{F}\left(\vartheta_{1}\right)-\mathfrak{F}\left(\vartheta_{2}\right)}{\cos\vartheta_{1}-\cos\vartheta_{2}}\,\mathrm{d}\vartheta_{1}\,\mathrm{d}\vartheta_{2}
		\]
		\[
		=\pi\left\{ 2\sum_{k=1}^{\infty}\alpha_{2k}\sum_{n=0}^{k-1}\dfrac{(-1)^{n}}{2n+1}+\dfrac{\pi}{4}\sum_{k=0}^{\infty}\alpha_{2k+1}+\pi\sum_{k=1}^{\infty}\alpha_{2k+1}\sum_{n=1}^{k}\dfrac{(-1)^{n+1}}{n}\sum_{m=0}^{n-1}\dfrac{1}{2m+1}\right\}. 
		\]
	\end{thm}
	As a consequence, we have the following corollary. 
	\begin{cor}
		Assume that $\mathfrak{F}$  has Fourier series expansion
		\[
		\mathfrak{F}\left(\vartheta\right)=\sum_{k=0}^{\infty}\alpha_{2k}\cos\left(2k\vartheta\right),
		\]
		then the following identity holds
		\[
		\int_{0}^{\pi/2}\int_{0}^{\pi/2}\dfrac{\mathfrak{F}\left(\vartheta_{1}\right)-\mathfrak{F}\left(\vartheta_{2}\right)}{\cos\vartheta_{1}-\cos\vartheta_{2}}\,\mathrm{d}\vartheta_{1}\,\mathrm{d}\vartheta_{2}=2\pi\sum_{k=1}^{\infty}\alpha_{2k}\sum_{n=0}^{k-1}\dfrac{\left(-1\right)^{n}}{2n+1}.
		\]
	\end{cor}
	An application of this corollary is provided by the first integral in Proposition \ref{cosh_int}, namely
	\[
	\int_{0}^{\pi/2}\int_{0}^{\pi/2}\dfrac{\log\left(\cos x\right)-\log\left(\cos z \right)}{\cos x-\cos z}\,\mathrm{d}x\,\mathrm{d}z =2\pi G.
	\]
	In this particular case, we have $\alpha_0=-\log 2,$ $\alpha_{2k}=(-1)^{k-1}/k$ 
	and $\alpha_{2k+1}=0,$ producing the expression of the Catalan constant $G$ as the alternate MZV of depth $2$,
	\begin{equation}
		\label{MZV11}    
		G=\sum_{k=1}^{\infty} \frac{\left(-1\right)^{k-1}}{k}
		\sum_{n=0}^{k-1}
		\frac{\left(-1\right)^{n}}{2n+1}.
	\end{equation}
	Two important remarks are worth mentioning at this point:
	\begin{enumerate}
		\item Representation \eqref{MZV11} can be obtained by expanding the integrand in the evaluation
		\[
		\int_{0}^{1} \frac{\log \left( 1-x^2 \right)}{1+x^2} \mathrm{d}x = \frac{\pi}{4} \log 2 -G
		\]
		that results from the addition of identities 4.291.8 and 4.291.10 in \cite{Gradshteyn}.
		\item There is an appealing parallel between identity \eqref{MZV11} rewritten as
		\[
		\sum_{k=0}^{\infty}\dfrac{\left(-1\right)^{k}}{\left(2k+1\right)^{2}}=\sum_{k= 1}^{\infty} \frac{\left(-1\right)^{k-1}}{k}
		\sum_{n=0}^{k-1}
		\frac{\left(-1\right)^{n}}{2n+1},
		\]
		and Euler's famous identity
		\[
		\sum_{k=1}^{\infty}\dfrac{1}{k^{3}}=\sum_{k = 1}^{\infty} \frac{1}{k^2}
		\sum_{n=1}^{k-1}
		\frac{1}{n}.
		\]
	\end{enumerate}
	
	Finally, 
	as an application of Theorem \ref{thm:generating_function}, the
	value of the original integral \eqref{double_integral} is recovered
	as follows.
	\begin{claim}
		\label{claim} The substitution $\alpha_{k}\mapsto\alpha^{k}$ in
		Theorem $\ref{thm:generating_function}$ allows us to recover Proposition
		$\ref{double_integral}$. 
	\end{claim}
	The proof of this claim is provided in Subsection \ref{subsec:Proof-of-Claim}.
	
	\section{\label{sec:A-Triple-Integral}A Triple Integral and a Triple Sum
		Version}
	A triple integral and a triple sum version of integral \eqref{double_integral}
	are now produced.
	\begin{thm}
		\label{thm:representations}The integral \eqref{double_integral} is
		represented as 
		\begin{align*}&\frac{1}{8}\int_{0}^{\pi/2}\int_{0}^{\pi/2}\dfrac{\log\left(\cos\left(\dfrac{x}{2}\right)\right)-\log\left(\cos\left(\dfrac{z}{2}\right)\right)}{\cos x-\cos z}\,\mathrm{d}x\,\mathrm{d}z \\&=\int_{0}^{1}\int_{0}^{z}\int_{0}^{y}\frac{\mathrm{d}x\,\mathrm{d}y\,\mathrm{d}z}{\left(1-x^{2}\right)\left(1+y^{2}\right)\left(1+z^{2}\right)}\\&=\sum_{\underset{k,l,m \mathrm{\,\, odd}}{0<k<l<m}}\frac{\left(-1\right)^{\frac{k+m}{2}-1}}{klm}.\end{align*}
	\end{thm}
	We notice that the nested triple series is an alternating multiple $t$-value
	as defined by Hoffman in \cite{Hoffman}: 
	\[
	t\left(k_{1},k_{2},\ldots,k_{n};z_{1},z_{2},\ldots,z_{n}\right)=\sum_{\underset{\mathrm{all\,odd}}{0<\omega_{1}<\cdots<\omega_{n}}
	}\dfrac{(-1)^{\omega_{1}z_{1}+\cdots+\omega_{n}z_{n}}}{\omega_{1}^{k_{1}}\omega_{2}^{k_{2}}\ldots\omega_{n}^{k_{n}}};
	\]
	In our case, this is
	\[
	-\,t\left(1,1,1;\frac{1}{2},0,\frac{1}{2}\right).
	\]
	As a direct consequence of this theorem, the following identity is
	obtained, for which an independent proof is produced in Section \ref{subsec:Proof-of-Proposition 8}.
	\begin{prop}
		\label{catalan_triple} The Catalan's constant is equal to 
		\[
		G=\dfrac{1}{\pi}\left\{ 8\int_{0}^{1}\int_{0}^{\vartheta_{3}}\int_{0}^{\vartheta_{2}}\dfrac{\mathrm{d}\vartheta_{1}\,\mathrm{d}\vartheta_{2}\,\mathrm{d}\vartheta_{3}}{\left(1-\vartheta_{1}^{2}\right)\left(1+\vartheta_{2}^{2}\right)\left(1+\vartheta_{3}^{2}\right)}+\dfrac{7}{4}\,\zeta(3)\right\}. 
		\]
	\end{prop}
	As a final remark in this section, it is noticed that the identity
	between the triple integral and the triple sum in Theorem \ref{thm:representations}
	can be obtained as a consequence of  simple computation
	\begin{align*}
		&       \int_{0}^{1}\int_{0}^{\vartheta_{3}}\int_{0}^{\vartheta_{2}}\dfrac{\mathrm{d}\vartheta_{1}\,\mathrm{d}\vartheta_{2}\,\mathrm{d}\vartheta_{3}}{\left(1-\vartheta_{1}^{2}\right)\left(1+\vartheta_{2}^{2}\right)\left(1+\vartheta_{3}^{2}\right)} \\&=\int_{0}^{1}\int_{0}^{\vartheta_{3}}\int_{0}^{\vartheta_{2}}\sum_{k=0}^{\infty}\vartheta_{1}^{2k}\sum_{m=0}^{\infty}(-1)^{m}\,\vartheta_{2}^{2m}\sum_{n=0}^{\infty}(-1)^{n}\,\vartheta_{3}^{2n}\,\mathrm{d}\vartheta_{1}\,\mathrm{d}\vartheta_{2}\,\mathrm{d}\vartheta_{3}\\
		& =\sum_{k=0}^{\infty}\sum_{m=0}^{\infty}\sum_{n=0}^{\infty}(-1)^{m+n}\int_{0}^{1}\int_{0}^{\vartheta_{3}}\int_{0}^{\vartheta_{2}}\vartheta_{1}^{2k}\,\vartheta_{2}^{2m}\,\vartheta_{3}^{2n}\,\mathrm{d}\vartheta_{1}\,\mathrm{d}\vartheta_{2}\,\mathrm{d}\vartheta_{3}\\
		& =\sum_{k=0}^{\infty}\sum_{m=0}^{\infty}\sum_{n=0}^{\infty}\dfrac{(-1)^{m+n}}{(2k+1)(2k+2n+2)(2k+2n+2m+3)}\\
		& =\sum_{\underset{\text{$k,l,m$  odd}}{0<k<l<m}}\frac{\left(-1\right)^{\frac{k+m}{2}-1}}{klm}.
	\end{align*}

	\section{\label{sec:Integration-over}Integration over $\left[0,\pi\right]^{2}$}
	
	Experimenting with variations on the original double integral
	\eqref{double_integral}, we changed the 
	domain of integration from $\left[0,\pi/2\right]^{2}$
	to $\left[0,\pi\right]^{2}$, considering integrals of the form
	\begin{equation}
		\int_{0}^{\pi}\int_{0}^{\pi}\frac{\mathfrak{F}\left(\cos\left(\vartheta_{1}\right)\right)-\mathfrak{F}\left(\cos\left(\vartheta_{2}\right)\right)}{\cos\vartheta_{1}-\cos\vartheta_{2}}\,\mathrm{d}\vartheta_{1}\,\mathrm{d}\vartheta_{2},\label{integral0Pi}
	\end{equation}
	and were rewarded with the following evaluations and their specializations.
	
	Here $\mathrm{Li}_{n}(z)$ represents the polylogarithm function defined by
	\[
	\mathrm{Li}_{n}\left(z\right)=\sum_{k=1}^{\infty}\frac{z^{k}}{k^{n}},
	\]
	and $\zeta(n)=\mathrm{Li}_{n}\left(1\right)$ denotes the Riemann
	zeta function.
	\begin{prop}
		\label{polylog_int} The following identity holds
		\[
		\int_{0}^{\pi}\int_{0}^{\pi}\dfrac{\mathrm{Li}_{n}\left(z\cos\vartheta_{1}\right)-\mathrm{Li}_{n}\left(z\cos\vartheta_{2}\right)}{\cos\vartheta_{1}-\cos\vartheta_{2}}\,\mathrm{d}\vartheta_{1}\,\mathrm{d}\vartheta_{2}=\dfrac{\pi^{2}}{2}\left[\mathrm{Li}_{n}\left(z\right)-\mathrm{Li}_{n}\left(-z\right)\right],
		\]
		and its specialization at $z=1$ produces 
		\[
		\int_{0}^{\pi}\int_{0}^{\pi}\dfrac{\mathrm{Li}_{n}\left(\cos\vartheta_{1}\right)-\mathrm{Li}_{n}\left(\cos\vartheta_{2}\right)}{\cos\vartheta_{1}-\cos\vartheta_{2}}\,\mathrm{d}\vartheta_{1}\,\mathrm{d}\vartheta_{2}=\pi^{2}\left(1-2^{-n}\right)\zeta(n),
		\]
		and, with $n=2$ and $z=\phi=(\sqrt{5}-1)/2$ and the golden ratio $\varphi=(1+\sqrt{5})/2$, 
		\[
		\int_{0}^{\pi}\int_{0}^{\pi}\dfrac{\mathrm{Li}_{2}\left(\phi\cos\vartheta_{1}\right)-\mathrm{Li}_{2}\left(\phi\cos\vartheta_{2}\right)}{\cos\vartheta_{1}-\cos\vartheta_{2}}\,\mathrm{d}\vartheta_{1}\,\mathrm{d}\vartheta_{2}=\dfrac{\pi^{4}}{12}-\dfrac{3}{2}\log^{2}\left(\varphi\right).
		\]
	\end{prop}
	In the case $\mathfrak{F}\left(z\right)=\log\left(1+z\right)$, 
	we obtain the following integral evaluation.
	\begin{prop} The following identity holds
		\label{lol_int} 
		\[
		\int_{0}^{\pi}\int_{0}^{\pi}\frac{\log\left(1+z\cos\vartheta_{1}\right)-\log\left(1+z\cos\vartheta_{2}\right)}{\cos\vartheta_{1}-\cos\vartheta_{2}}\,\mathrm{d}\vartheta_{1}\,\mathrm{d}\vartheta_{2}=\frac{\pi^{2}}{2}\log\left(\frac{1+z}{1-z}\right)
		.\]
	\end{prop}
	
	But soon it was realized that the double integration \eqref{integral0Pi}
	over $\left[0,\pi\right]^{2}$ simply extracts the odd and even part
	of the function $\mathfrak{F}$, as shown in the following theorem. 
	\begin{thm}
		\label{analytic_thm} If $\mathfrak{F}\left(z\right)=\sum_{n\geqslant0}a_{n}z^{n}$
		is an analytic function, then we have
		\begin{align*}
			\frac{1}{\pi^{2}}\int_{0}^{\pi}\int_{0}^{\pi}\frac{\mathfrak{F}\left(\cos\vartheta_{1}\right)-\mathfrak{F}\left(\cos\vartheta_{2}\right)}{\cos\vartheta_{1}-\cos\vartheta_{2}}\,\mathrm{d}\vartheta_{1}\,\mathrm{d}\vartheta_{2}=\frac{\mathfrak{F}\left(z\right)-\mathfrak{F}\left(-z\right)}{2}=\sum_{n=0}^{\infty}a_{2n+1}z^{2n+1},
		\end{align*}
		and 
		\begin{align*}
			\frac{1}{\pi^{2}}\int_{0}^{\pi}\int_{0}^{\pi}\frac{\cos\vartheta_{1}\mathfrak{F}\left(\cos\vartheta_{1}\right)-\cos\vartheta_{2}\mathfrak{F}\left(\cos\vartheta_{2}\right)}{\cos\vartheta_{1}-\cos\vartheta_{2}}\mathrm{d}\vartheta_{1}\mathrm{d}\vartheta_{2}=\frac{\mathfrak{F}\left(z\right)+\mathfrak{F}\left(-z\right)}{2}=\sum_{n=0}^{\infty}a_{2n}z^{2n}.
		\end{align*}
	\end{thm}
	
	This result suggests that, due to
	the symmetry of the integration domain, versions over $\left[0,\pi\right]^{2}$
	of the previous integrals are somewhat less interesting than expected.
	
	\section{\label{sec:A-parameterized-version}A Parameterized Version}
	
	It sometimes happens that an integral is easier to compute when considered
	as a specialization of a parameterized family of integrals. This principle
	is illustrated next in the case of integral \eqref{double_integral}
	with the following result, producing another proof of Proposition \ref{Prop1}.
	\begin{prop}
		\label{one_parameter} The following identity holds
		\[
		\int_{0}^{\pi/2}\int_{0}^{\pi/2}\frac{\log\left(1+\cos z\cos\vartheta_{1}\right)-\log\left(1+\cos z\cos\vartheta_{2}\right)}{\cos\vartheta_{1}-\cos\vartheta_{2}}\,\mathrm{d}\vartheta_{1}\,\mathrm{d}\vartheta_{2}\]
		\begin{equation*}
			=2\imath z\left[\Li_{2}\left(\imath e^{\imath z}\right)-\Li_{2}\left(-\imath e^{\imath z}\right)\right]+2\left[\Li_{3}\left(-\imath e^{\imath z}\right)-\Li_{3}\left(\imath e^{\imath z}\right)\right]-z^{2}\log\left(\dfrac{1-e^{\imath z}}{1+e^{\imath z}}\right)+2\pi G.
		\end{equation*}
		The specialization at $z=0$ produces 
		\[
		\int_{0}^{\pi/2}\int_{0}^{\pi/2}\frac{\log\left(1+\cos\vartheta_{1}\right)-\log\left(1+\cos\vartheta_{2}\right)}{\cos\vartheta_{1}-\cos\vartheta_{2}}\,\mathrm{d}\vartheta_{1}\,\mathrm{d}\vartheta_{2}=2\pi G-\dfrac{7}{2}\,\zeta(3),
		\]
		which is equivalent to equation \eqref{double_integral}.
	\end{prop}
	The proof of this Proposition is given in Subsection \ref{one_parameter_proof}.
	
	\section{\label{sec:A-simple-integral}A Simple Integral Representation}
	
	In an attempt to generalize integral \eqref{double_integral}, let us
	now consider the case of two functions $\mathfrak{F}$ and $H$ and the 
	integral 
	\[
	\int_{\alpha}^{\beta}\int_{\alpha}^{\beta}\dfrac{\mathfrak{F}\left(zH\left(\vartheta_{1}\right)\right)-\mathfrak{F}\left(zH\left(\vartheta_{2}\right)\right)}{H\left(\vartheta_{1}\right)-H\left(\vartheta_{2}\right)}\,\mathrm{d}\vartheta_{1}\,\mathrm{d}\vartheta_{2}.
	\]
	It appears that, using the notion of Hadamard product, this integral can be expressed
	as a single integral as follows. 
	\begin{thm}
		\label{lil_general} Assume that $\mathfrak{F}\left(z\right)$ is an analytic function in the unit disk
		then 
		\begin{align*}&
			\int_{0}^{\pi/2}\int_{0}^{\pi/2}\dfrac{\mathfrak{F}\left(z\cos\vartheta_{1}\right)-\mathfrak{F}\left(z\cos\vartheta_{2}\right)}{\cos\vartheta_{1}-\cos\vartheta_{2}}\,\mathrm{d}\vartheta_{1}\,\mathrm{d}\vartheta_{2}
			\\&
			=\dfrac{1}{2\pi}\int_{0}^{2\pi}\mathfrak{F}\left(e^{\imath\vartheta}\sqrt{z}\right)\dfrac{e^{-\imath\vartheta}\sqrt{z}}{1-ze^{-2\imath\vartheta}}\arccos^{2}\left(\sqrt{1-ze^{-2\imath\vartheta}}\right)\mathrm{d}\vartheta.
		\end{align*}
		Moreover, assuming that the limit $z \to 1$ exists,
		\[
		\int_{0}^{\pi/2}\int_{0}^{\pi/2}\dfrac{\mathfrak{F}\left(\cos\vartheta_{1}\right)-\mathfrak{F}\left(\cos\vartheta_{2}\right)}{\cos\vartheta_{1}-\cos\vartheta_{2}}\,\mathrm{d}\vartheta_{1}\,\mathrm{d}\vartheta_{2}=\dfrac{1}{4\imath\pi}\int_{0}^{2\pi}\dfrac{\mathfrak{F}\left(e^{\imath\vartheta}\right)}{\sin\left(\vartheta\right)}\arccos^{2}\left(\sqrt{1-e^{-2\imath\vartheta}}\right)\mathrm{d}\vartheta.
		\]
	\end{thm}
	
	Theorem \ref{lil_general} can be easily extended to the following result.
	\begin{thm}
		\label{general_single} Given two functions $\mathfrak{F}\left(z\right)$ and $H\left(z\right)$, defining
		the generating function 
		\[
		G\left(z\right)=\left(\int_{\alpha}^{\beta}\dfrac{1}{1-zH\left(\vartheta\right)}\,\mathrm{d}\vartheta\right)^{2},
		\]
		and assuming that 
		\[
		G\left(z\right)=\sum_{n=0}^{\infty}\mathfrak{B}_{n}z^{n},\quad \mathfrak{F}\left(z\right)=\sum_{n=0}^{\infty}\mathfrak{R}_{n}z^{n},
		\]
		then the following identity holds
		\[
		\int_{\alpha}^{\beta}\int_{\alpha}^{\beta}\dfrac{\mathfrak{F}\left(zH\left(\vartheta_{1}\right)\right)-\mathfrak{F}\left(zH\left(\vartheta_{2}\right)\right)}{H\left(\vartheta_{1}\right)-H\left(\vartheta_{2}\right)}\,\mathrm{d}\vartheta_{1}\,\mathrm{d}\vartheta_{2}=\sum_{n=1}^{\infty}\mathfrak{R}_{n}\mathfrak{B}_{n-1}z^{n},
		\]
		so that the general double integral can be represented as a simple integral:
		\[
		\int_{\alpha}^{\beta}\int_{\alpha}^{\beta}\dfrac{\mathfrak{F}\left(zH\left(\vartheta_{1}\right)\right)-\mathfrak{F}\left(zH\left(\vartheta_{2}\right)\right)}{H\left(\vartheta_{1}\right)-H\left(\vartheta_{2}\right)}\,\mathrm{d}\vartheta_{1}\,\mathrm{d}\vartheta_{2}=\dfrac{1}{2\pi}\int_{0}^{2\pi}\mathfrak{F}\left(e^{\imath\vartheta}\sqrt{z}\right)e^{-\imath\vartheta}\sqrt{z}\,G\left(e^{\imath\vartheta}\sqrt{z}\right)\mathrm{d}\vartheta.
		\]
	\end{thm}

	\section{\label{sec:multivariate-extensions}Multivariate Extensions and a
		link with The Stieltjes Transforms}
	One may wonder how integral \eqref{double_integral} can be extended
	to more than two variables. One possibility is given by the following
	result that establishes a link with Stieltjes transforms. 
	
	Consider, for $1\leqslant i\leqslant n+1,$ a set of functions 
	\[
	\mathfrak{f}_{i}:\mathbb{R}\to\mathbb{R},
	\]
	and assume the existence on a suitable domain of their transforms 
	\[
	\hat{\mathfrak{f}}_{i}\left(z\right)=\int_{\mathbb{R}}\frac{\mathfrak{f}_{i}\left(x\right)}{1+xz}\,\mathrm{d}x.
	\]
	Notice that $\hat{\mathfrak{f}}_{i}$ is a scaled version the Stieltjes transform
	$S_{\mathfrak{f}_{i}}$ of $\mathfrak{f}_{i}$ defined as
	\[
	S_{\mathfrak{f}_{i}}\left(z\right)=\int_{\mathbb{R}}\frac{\mathfrak{f}_{i}\left(x\right)}{z-x}\,\mathrm{d}x,
	\]
	since 
	\[
	\hat{\mathfrak{f}}_{i}\left(z\right)=-\frac{1}{z}\,S_{\mathfrak{f}_{i}}\left(-\frac{1}{z}\right).
	\]
	The integrals of the functions $\mathfrak{f}_{i}$ and of their transforms $\hat{\mathfrak{f}}_{i}$ are related
	through the following duality identity, which is one the main results of this study. 
	\begin{thm}
		\label{thm:Stieltjes}Define the function
		\[
		\hat{\mathfrak{F}}_{n+1}\left(x_{1},x_2,\ldots,x_{n}\right)=\sum_{i=1}^{n}\frac{x_{i}^{n-1}}{\prod_{j\ne i}\left(x_{i}-x_{j}\right)}\,\hat{\mathfrak{f}}_{n+1}\left(x_{i}\right),
		\]
		then the following duality identity holds between the functions $\mathfrak{f}_{i}$
		and their transforms $\hat{\mathfrak{f}}_{i}$
		\[
		\int_{\mathbb{R}^{n}}\hat{\mathfrak{F}}_{n+1}\left(x_{1},x_2,\ldots,x_{n}\right)\prod_{i=1}^{n}\mathfrak{f}_{i}\left(x_{i}\right)\mathrm{d}x_{i}=\int_{\mathbb{R}}\mathfrak{f}_{n+1}\left(x_{n+1}\right)\prod_{i=1}^{n}\hat{\mathfrak{f}}_{i}\left(x_{n+1}\right)\mathrm{d}x_{n+1}.
		\]
	\end{thm}
	
	Note that this identity relates an $n$-variate integral to a single-variable
	integral. However, observe that each term $\hat{\mathfrak{f}}_{i}$ is itself an integral.
	
	Two specializations of this theorem are shown next.
	\begin{prop}
		\label{cos_multivariant} For $n\geqslant 2,$ we have 
		\[
		\int_{\left[0,\pi/2\right]^{n}}\sum_{i=1}^{n}\frac{\cos^{n-2}\left(\vartheta_{i}\right)\log\cos\left(\vartheta_{i}/2\right)}{\prod_{j\ne i}\left(\cos\vartheta_{i}-\cos\vartheta_{j}\right)}\prod_{i=1}^{n}\mathrm{d}\vartheta_{i}=\int_{0}^{\pi/2}\dfrac{z^{n}}{\sin^{n-1}z}\,\mathrm{d}z.
		\]
		The special case $n=2$ recovers the double integral \eqref{double_integral}. 
		
		The case $n=3$ reads
		\begin{align*}
			\int_{0}^{\pi/2}\int_{0}^{\pi/2}\int_{0}^{\pi/2}&\frac{\cos\vartheta_{1}\log\left(\cos\left(\vartheta_{1}/2\right)\right)}{\left(\cos\vartheta_{1}-\cos\vartheta_{2}\right)\left(\cos\vartheta_{1}-\cos\vartheta_{3}\right)}+\frac{\cos\vartheta_{2}\log\left(\cos\left(\vartheta_{2}/2\right)\right)}{\left(\cos\vartheta_{2}-\cos\vartheta_{1}\right)\left(\cos\vartheta_{2}-\cos\vartheta_{3}\right)}\\
			+\,&\frac{\cos\vartheta_{3}\log\left(\cos\left(\vartheta_{3}/2\right)\right)}{\left(\cos\vartheta_{3}-\cos\vartheta_{1}\right)\left(\cos\vartheta_{3}-\cos\vartheta_{2}\right)}\,\mathrm{d}\vartheta_{1}\,\mathrm{d}\vartheta_{2}\,\mathrm{d}\vartheta_{3}=\frac{3\pi^{2}}{8}\log4-\frac{21}{8}\,\zeta(3).
		\end{align*}
	\end{prop}
	
	It is noticed that an explicit value of the integral $$J_{n}=\int_{0}^{\pi/2}\dfrac{z^{n}}{\sin^{n-1}z}\,\mathrm{d}z,$$
	for an arbitrary value of the parameter $n$ does not seem to be available. However, using integration by parts, first few values of $J_n$ are computed as
	\[
	J_{1}=\frac{\pi^{2}}{8},\quad J_{2}=2\pi G-\frac{7}{2}\,\zeta(3),\quad J_{3}=\frac{3\pi^{2}}{8}\log4-\frac{21}{8}
	\,\zeta(3),
	\]
	\[
	J_{4}=\frac{3}{4}\left(31
	\,\zeta(5)-28\,\zeta(3)\right)+\frac{\pi}{2}\,G\left(\pi^{2}+24\right)-\frac{\pi^{2}}{4}+\frac{1}{128}\left(\psi^{\left(3\right)}\left(\frac{3}{4}\right)-\psi^{\left(3\right)}\left(\frac{1}{4}\right)\right),
	\]
	where $\psi^{\left(3\right)}\left(z\right)$ represents the pentagamma function. 
	
	\begin{prop}
		\label{tan_multivariant} For $n\geqslant2,$ the following identity holds
		\[
		\int_{\left[0,\pi/2\right]^{n}}\sum_{i=1}^{n}\frac{\tan^{n-2}\vartheta_{i}\log\left(1+\tan\vartheta_{i}\right)}{\prod_{j\ne i}\left(\tan\vartheta_{i}-\tan\vartheta_{j}\right)}\prod_{i=1}^{n}\mathrm{d}\vartheta_{i}=\int_{0}^{1}\frac{\left(\pi/2+z\log z\right)^{n}}{\left(1+z^{2}\right)^{n}}\,\mathrm{d}z.
		\]
		The specialization $n=2$ produces
		\[
		\int_{0}^{\pi/2}\int_{0}^{\pi/2}\frac{\log\left(1+\tan\vartheta_{1}\right)-\log\left(1+\tan\vartheta_{1}\right)}{\tan\vartheta_{1}-\tan\vartheta_{2}}\mathrm{d}\vartheta_{1}\mathrm{d}\vartheta_{2}=\frac{\pi}{16}\left(\pi+\pi^{2}-4\log2\right)-G.
		\]
	\end{prop}

	\section{\label{proof} Proofs}
	
	\subsection{\label{subsec:Proof Cosh}Proof of proposition \ref{cosh_int}}
	The first integral
	\[
	I=\int_{0}^{\pi/2}\int_{0}^{\pi/2}\frac{\log\cos x-\log\cos z}{\cos x-\cos z}\,\mathrm{d}x\,\mathrm{d}z,
	\]
	is first transformed as
	\[
	I=\frac{1}{2}\int_{0}^{\pi/2}\int_{0}^{\pi/2}\frac{\log\left(1-\sin^{2}x\right)-\log\left(1-\sin^{2}z\right)}{\cos x-\cos z}\,\mathrm{d}x\,\mathrm{d}z,
	\]
	and rewritten as the triple integral
	\[
	I=\frac{1}{2}\int_{0}^{1}\int_{0}^{\pi/2}\int_{0}^{\pi/2}\dfrac{\frac{\sin^{2}z}{1-t\sin^{2}z}-\frac{\sin^{2}x}{1-t\sin^{2}x}}{\cos x-\cos z}\,\mathrm{d}x\,\mathrm{d}z\,\mathrm{d}t,
	\]
	which simplifies as
	\begin{align*}
		\int_{0}^{\pi/2}\int_{0}^{\pi/2}\frac{\log\cos x-\log\cos z}{\cos x-\cos z}\,\mathrm{d}x\,\mathrm{d}z & =\frac{1}{2}\int_{0}^{1}\int_{0}^{\pi/2}\int_{0}^{\pi/2}\frac{\frac{\sin^{2}x-\sin^{2}z}{\left(1-t\sin^{2}x\right)\left(1-t\sin^{2}z\right)}}{\cos x-\cos z}\,\mathrm{d}x\,\mathrm{d}z\,\mathrm{d}t\\
		&=\frac{1}{2}\int_{0}^{1}\int_{0}^{\pi/2}\int_{0}^{\pi/2}\frac{\cos x+\cos z}{\left(1-t\sin^{2}x\right)\left(1-t\sin^{2}z\right)}\,\mathrm{d}x\,\mathrm{d}z\,\mathrm{d}t,
	\end{align*}
	and, by symmetry, as
	\[
	I=\int_{0}^{1}\int_{0}^{\pi/2}\int_{0}^{\pi/2}\frac{\cos x}{\left(1-t\sin^{2}x\right)\left(1-t\sin^{2}z\right)}\,\mathrm{d}x\,\mathrm{d}z\,\mathrm{d}t.
	\]
	Integrating over $z$ and then over $x$ produces
	\begin{align*}
		I  =\int_{0}^{1}\int_{0}^{\pi/2}\frac{\cos x}{\left(1-t\sin^{2}x\right)}\frac{\pi}{2\sqrt{1-t}}\,\mathrm{d}x\,\mathrm{d}t =\frac{\pi}{4}\int_{0}^{1}\frac{1}{\sqrt{t\left(1-t\right)}}\log\left(\frac{1+\sqrt{t}}{1-\sqrt{t}}\right)\mathrm{d}t.
	\end{align*}
	The change of variables $t=\cos^{2}\left(\dfrac{u}{2}\right)$
	and the representation $G=\int_{0}^{\frac{\pi}{4}}\log\cot u\,\mathrm{d}u$
	produces
	\[
	I=\frac{\pi}{4}\int_{0}^{\pi}\log\left(\frac{1+\cos\left(\dfrac{u}{2}\right)}{1-\cos\left(\dfrac{u}{2}\right)}\right)\mathrm{d}u=\frac{\pi}{4}\int_{0}^{\pi}2\log\cot\left(\dfrac{u}{4}\right)\mathrm{d}u=2\pi G,
	\]
	which is the desired result.
	
	For the second integral, start with the identity
	\[
	\int_{1}^{\infty}\frac{1}{1+xz}\frac{\mathrm{d}x}{\sqrt{x^{2}-1}}=\frac{\mathrm{arcsec}\,z}{\sqrt{z^{2}-1}},\quad\left(0\leqslant z\leqslant1\right)
	,\]
	so that
	\[
	\int_{0}^{\infty}\int_{0}^{\infty}\frac{1}{1+x_{1}z}\frac{1}{1+x_{2}z}\frac{\mathrm{d}x_{1}}{\sqrt{x_{1}^{2}-1}}\frac{\mathrm{d}x_{2}}{\sqrt{x_{2}^{2}-1}}=\left(\frac{\mathrm{arcsec}\,z}{\sqrt{z^{2}-1}}\right)^{2},
	\]
	and
	\[
	\int_{0}^{1}\mathrm{d}z\int_{0}^{\infty}\int_{0}^{\infty}\frac{1}{1+x_{1}z}\frac{1}{1+x_{2}z}\frac{\mathrm{d}x_{1}}{\sqrt{x_{1}^{2}-1}}\frac{\mathrm{d}x_{2}}{\sqrt{x_{2}^{2}-1}}=\int_{0}^{1}\left(\frac{\mathrm{arcsec}\,z}{\sqrt{z^{2}-1}}\right)^{2}\mathrm{d}z.
	\]
	Since 
	\[
	\int_{0}^{1}\left(\frac{\mathrm{arcsec}\,z}{\sqrt{z^{2}-1}}\right)^{2}\mathrm{d}z=\dfrac{7}{2}\,\zeta(3),
	\]
	the change of variable $x_{1}=\cosh\vartheta_{1}$ and $x_{2}=\cosh\vartheta_{2}$
	produces
	\begin{align*}
		&\int_{0}^{1}\mathrm{d}z\int_{0}^{\infty}\int_{0}^{\infty}\frac{1}{1+x_{1}z}\frac{1}{1+x_{2}z}\frac{\mathrm{d}x_{1}}{\sqrt{x_{1}^{2}-1}}\frac{\mathrm{d}x_{2}}{\sqrt{x_{2}^{2}-1}}
		\\&	 =\int_{1}^{\infty}\int_{1}^{\infty}\frac{\log\left(1+x_{1}\right)-\log\left(1+x_{2}\right)}{x_{1}-x_{2}}\frac{\mathrm{d}x_{1}}{\sqrt{x_{1}^{2}-1}}\frac{\mathrm{d}x_{2}}{\sqrt{x_{2}^{2}-1}}
		\\& =\int_{0}^{\infty}\int_{0}^{\infty}\frac{\log\left(1+\cosh\vartheta_{1}\right)-\log\left(1+\cosh\vartheta_{2}\right)}{\cosh\vartheta_{1}-\cosh\vartheta_{2}}\,\mathrm{d}\vartheta_{1}\,\mathrm{d}\vartheta_{2}=\dfrac{7}{2}\,\zeta(3),
	\end{align*}
	which is the desired result. \QED
	\subsection{Proof of Proposition \ref{prop:Prop2}}
	\label{ProofProp2}
	Start from the identity 
	\[
	\log\left(\cos\left(\dfrac{\vartheta}{2}\right)\right)=\sum_{k=1}^{\infty}\dfrac{(-1)^{k-1}}{k}\cos\left(k\vartheta\right)-\log2,
	\]
	so that the double integral, which we will call $\mathfrak{D}$, is
	equal to 
	\[
	\mathfrak{D}=\int_{0}^{\pi/2}\int_{0}^{\pi/2}\dfrac{\log\left(\cos\left(\dfrac{x}{2}\right)\right)-\log\left(\cos\left(\dfrac{z}{2}\right)\right)}{\cos\left(x\right)-\cos\left(z\right)}\,\mathrm{d}x\,\mathrm{d}z
	\]
	
	\[
	\,\,\,\,\,\,=\sum_{k=1}^{\infty}\dfrac{(-1)^{k-1}}{k}\int_{0}^{\pi/2}\int_{0}^{\pi/2}\dfrac{\cos\left(kx\right)-\cos\left(kz\right)}{\cos x-\cos z}\,\mathrm{d}x\,\mathrm{d}z.
	\]
	Therefore we need to compute the sequence of integrals
	\[
	\mathfrak{T}_{k}=\int_{0}^{\pi/2}\int_{0}^{\pi/2}\dfrac{\cos\left(kx\right)-\cos\left(kz\right)}{\cos x-\cos z}\,\mathrm{d}x\,\mathrm{d}z,\quad\left(k\geqslant1\right),
	\]
	to obtain
	\[
	\mathfrak{D}=\sum_{k=1}^{\infty}\dfrac{(-1)^{k-1}}{k}\mathfrak{T}_{k}.
	\]

	\subsection*{Computation of $\mathfrak{T}_{2k}$}
	To start with, observe that 
	\[
	\mathfrak{T}_{2k}=\int_{0}^{\pi/2}\int_{0}^{\pi/2}\dfrac{\cos\left(2kx\right)-\cos\left(2kz\right)}{\cos x-\cos z}\,\mathrm{d}x\,\mathrm{d}z=\int_{0}^{\pi/2}\int_{0}^{\pi}\dfrac{\cos\left(2kx\right)-\cos\left(2kz\right)}{\cos x-\cos z}\,\mathrm{d}x\,\mathrm{d}z,
	\]
	and the fact that 
	\[
	\int_{-1}^{+1}\dfrac{T_{2k}\left(z_{1}\right)}{z_{1}-z_{2}}\dfrac{\mathrm{d}z_{1}}{\sqrt{1-z_{1}^{2}}}=\pi\mathfrak{Z}_{2k-1}\left(z_{2}\right),
	\]
	so that we have 
	\[
	\mathfrak{T}_{2k}=\pi\int_{0}^{\pi/2}\mathfrak{Z}_{2k-1}\left(\cos z\right)\mathrm{d}z=\pi\int_{0}^{\pi/2}\dfrac{\sin\left(2kz\right)}{\sin\left(z\right)}\,\mathrm{d}z=2\pi\sum_{\ell=1}^{k-1}\dfrac{\left(-1\right)^{\ell}}{2\ell+1},
	\]
	the last identity being given in \cite[Identity 2.5.12.7]{Xu}. Therefore we deduce
	that 
	\begin{equation}
		\mathfrak{T}_{2k}=\int_{0}^{\pi/2}\int_{0}^{\pi/2}\dfrac{\cos\left(2kx\right)-\cos\left(2kz\right)}{\cos x-\cos z}\,\mathrm{d}x\,\mathrm{d}z=2\pi\sum_{\ell=1}^{k-1}\dfrac{\left(-1\right)^{\ell}}{2\ell+1}.\label{prud_identity}
	\end{equation}

	\subsection*{Computation of $\mathfrak{T}_{2k+1}$}
	
	We can easily check that 
	\begin{align*}
		\frac{\cos\left(2n+1\right)x-\cos\left(2n+1\right)z}{\cos x-\cos z}-\frac{\cos\left(2n-1\right)x-\cos\left(2n-1\right)z}{\cos x-\cos z}=\sum_{\vert p\vert+\vert q\vert=2n}\cos\left(px+qz\right).
	\end{align*}
	Integration produces 
	\begin{align*}
		\int_{0}^{\pi/2}\int_{0}^{\pi/2}\cos\left(px+qz\right)\mathrm{d}x\,\mathrm{d}z & =\frac{4}{pq}\sin\left(\frac{p\pi}{4}\right)\sin\left(\frac{q\pi}{4}\right)\cos\left(\frac{p\pi+q\pi}{4}\right)\\
		& =\frac{\pi^{2}}{4}\,\mathrm{sinc}\left(\frac{p\pi}{4}\right)\mathrm{sinc}\left(\frac{q\pi}{4}\right)\cos\left(\frac{p\pi+q\pi}{4}\right).
	\end{align*}
	Therefore, we have
	\[
	\mathfrak{T}_{2n+1}=\mathfrak{T}_{2n-1}+\frac{\pi^{2}}{4}\sum_{\vert p\vert+\vert q\vert=2n}\mathrm{sinc}\left(\frac{p\pi}{4}\right)\mathrm{sinc}\left(\frac{q\pi}{4}\right)\cos\left(\frac{p\pi+q\pi}{4}\right),
	\]
	and with $\mathfrak{T}_{1}=\left(\pi/2\right)^{2}$, we deduce that
	\[
	\mathfrak{T}_{2n+1}=\left(\frac{\pi}{2}\right)^{2}\left[1+\sum_{k=1}^{n}\sum_{\vert p\vert+\vert q\vert=2k}\mathrm{sinc}\left(\frac{p\pi}{4}\right)\mathrm{sinc}\left(\frac{q\pi}{4}\right)\cos\left(\frac{p\pi+q\pi}{4}\right)\right].
	\]
	Let us first compute the sum 
	\[
	\sum_{p+q=2n}\mathrm{sinc}\left(\frac{p\pi}{4}\right)\mathrm{sinc}\left(\frac{q\pi}{4}\right)\cos\left(\frac{p\pi+q\pi}{4}\right)=\cos\left(\frac{n\pi}{2}\right)\sum_{p+q=2n}\mathrm{sinc}\left(\frac{p\pi}{4}\right)\mathrm{sinc}\left(\frac{q\pi}{4}\right)
	\]
	
	\[
	=\cos\left(\frac{n\pi}{2}\right)\left[2\,\mathrm{sinc}\left(\frac{n\pi}{2}\right)+\sum_{p=1}^{2n-1}\mathrm{sinc}\left(\frac{p\pi}{4}\right)\mathrm{sinc}\left(\frac{2n\pi-p\pi}{4}\right)\right].
	\]
	Next, we compute the sum 
	\[
	\sum_{p=1}^{2n-1}\dfrac{1}{p\left(2n-1\right)}\,\mathrm{sinc}\left(\frac{p\pi}{4}\right)\mathrm{sinc}\left(\frac{2n\pi-p\pi}{4}\right),
	\]
	linearization produces 
	\[
	\mathrm{sinc}\left(\frac{p\pi}{4}\right)\mathrm{sinc}\left(\frac{2n\pi-p\pi}{4}\right)=\dfrac{1}{2}\cos\left(\dfrac{n\pi-p\pi}{2}\right)-\dfrac{1}{2}\cos\left(\dfrac{n\pi}{2}\right),
	\]
	so that we have 
	\begin{align*}&
		\sum_{p=1}^{2n-1}\dfrac{1}{p\left(2n-1\right)}\,\mathrm{sinc}\left(\frac{p\pi}{4}\right)\mathrm{sinc}\left(\frac{2n\pi-p\pi}{4}\right)
		\\&
		=\dfrac{1}{2}\sum_{p=1}^{2n-1}\dfrac{1}{p\left(2n-1\right)}\cos\left(\dfrac{n\pi-p\pi}{2}\right)-\dfrac{1}{2}\cos\left(\dfrac{n\pi}{2}\right)\sum_{p=1}^{2n-1}\dfrac{1}{p\left(2n-1\right)}.
	\end{align*}
	A long but simple computation produces 
	\[
	\mathfrak{T}_{2n+1}-\mathfrak{T}_{2n-1}=\dfrac{\pi^{2}}{4}\sum_{p+q=2n}\mathrm{sinc}\left(\frac{p\pi}{4}\right)\mathrm{sinc}\left(\frac{q\pi}{4}\right)\cos\left(\frac{p\pi+q\pi}{4}\right)
	\]
	\[
	=\frac{2\left(-1\right)^{n+1}}{n}\left[2H_{2n-1}-H_{n-1}\right]=\frac{4\left(-1\right)^{n+1}}{n}\sum_{m=0}^{n-1}\frac{1}{2m+1},
	\]
	so that we have
	\[
	\mathfrak{T}_{2n+1}=\frac{\pi^{2}}{4}+4\sum_{k=1}^{n}\left(-1\right)^{k+1}\frac{H_{2k-1}}{k}-2\sum_{k=1}^{n}\left(-1\right)^{k+1}\frac{H_{k-1}}{k}.
	\]
	Putting all things together produces
	\begin{equation}
		\int_{0}^{\pi/2}\int_{0}^{\pi/2}\dfrac{\cos\left(\left(2k+1\right)x\right)-\cos\left(\left(2k+1\right)z\right)}{\cos x-\cos z}\,\mathrm{d}x\,\mathrm{d}z=\frac{\pi^{2}}{4}+\sum_{k=1}^{n}\frac{4\left(-1\right)^{k+1}}{k}\sum_{m=0}^{n-1}\frac{1}{2m+1}\label{eq:e2n+1}.
	\end{equation}
	Finally, we need to compute the sum 
	\[
	\mathfrak{D}=\sum_{k=1}^{\infty}\dfrac{(-1)^{k-1}}{k}\int_{0}^{\pi/2}\int_{0}^{\pi/2}\dfrac{\cos\left(kx\right)-\cos\left(kz\right)}{\cos x-\cos z}\,\mathrm{d}x\,\mathrm{d}z.
	\]
	Since we already know that 
	\begin{equation}
		\mathfrak{T}_{2k}=2\pi\sum_{\ell=1}^{k-1}\dfrac{\left(-1\right)^{\ell}}{2\ell+1},\quad\mathfrak{T}_{2k+1}=\frac{\pi^{2}}{4}+4\sum_{k=1}^{n}\frac{\left(-1\right)^{k+1}}{k}\sum_{m=0}^{k-1}\frac{1}{2m+1},\label{computations_integral-1}
	\end{equation}
	we thus need to evaluate 
	\[
	\mathfrak{C}=\frac{\pi^{2}}{4}-\sum_{k=1}^{\infty}\left\{ \frac{\pi}{k}\sum_{\ell=1}^{k-1}\dfrac{\left(-1\right)^{\ell}}{2\ell+1}+\frac{\pi^{2}}{4\left(2k+1\right)}+\frac{4}{2k+1}\sum_{\ell=1}^{k}\frac{\left(-1\right)^{\ell+1}}{\ell}\sum_{m=0}^{\ell-1}\frac{1}{2m+1}\right\}. 
	\]
	Note that the terms can not be computed separately since for example,
	the sum 
	\[
	\sum_{k=1}^{\infty}\dfrac{1}{2k+1},
	\]
	is not convergent. However, notice that 
	\[
	\sum_{\ell=1}^{k-1}\dfrac{\left(-1\right)^{\ell}}{2\ell+1}=\frac{\pi}{4}-\sum_{l\geqslant k}\dfrac{\left(-1\right)^{\ell}}{2\ell+1},
	\]
	and that 
	\[
	4\sum_{\ell=1}^{k}\frac{\left(-1\right)^{\ell+1}}{\ell}\sum_{m=0}^{\ell-1}\frac{1}{2\ell+1}=\frac{\pi^{2}}{16}-4\sum_{\ell\geqslant k+1}\frac{\left(-1\right)^{\ell+1}}{\ell}\sum_{m=0}^{\ell-1}\frac{1}{2m+1},
	\]
	which yields 
	\begin{align*}
		\mathfrak{D} & =\frac{\pi^{2}}{4}-\sum_{k=1}^{\infty}\left[\dfrac{\pi}{k}\left(\frac{\pi}{4}-\sum_{\ell\geqslant k}\frac{\left(-1\right)^{\ell}}{2\ell+1}\right)+\frac{\pi^{2}}{4\left(2k+1\right)}\right]\\
		& +\sum_{k=1}^{\infty}\left[\dfrac{4}{2k+1}\left(\frac{\pi^{2}}{16}-4\sum_{\ell\geqslant k+1}\frac{\left(-1\right)^{\ell+1}}{\ell}\sum_{m=0}^{\ell-1}\frac{1}{2m+1}\right)\right]\\
		& =\frac{\pi^{2}}{4}\left[1+\sum_{k=1}^{\infty}\left(\frac{1}{2k+1}+\frac{1}{2k+1}-\frac{1}{k}\right)\right]\\
		& +\sum_{k=1}^{\infty}\frac{\pi}{k}\sum_{\ell\geqslant k}\frac{\left(-1\right)^{\ell}}{2\ell+1}-\sum_{k=1}^{\infty}\frac{4}{2k+1}\sum_{\ell\geqslant k+1}\frac{\left(-1\right)^{\ell+1}}{\ell}\sum_{m=0}^{\ell-1}\frac{1}{2m+1}.
	\end{align*}
	The first sum is easily computed as 
	\[
	\sum_{k=1}^{\infty}\left(\frac{1}{2k+1}+\frac{1}{2k+1}-\frac{1}{k}\right)=\sum_{k=1}^{\infty}\left(\frac{2}{2k+1}-\frac{1}{k}\right)=2\log\left(2\right)-2,
	\]
	while the second sum is expressed as the integral 
	\[
	\sum_{k=1}^{\infty}\frac{\pi}{k}\sum_{\ell\geqslant k}\frac{\left(-1\right)^{\ell}}{2\ell+1}=\pi\sum_{k=1}^{\infty}\frac{1}{k}\sum_{\ell=0}^{\infty}\frac{\left(-1\right)^{\ell+k}}{2\ell+2k+1}=\pi\sum_{k=1}^{\infty}\frac{\left(-1\right)^{k}}{k}\sum_{\ell=0}^{\infty}\frac{1}{2\ell+2k+1}
	\]
	\[
	\,\,\,\,\,\,=-\,\pi\int_{0}^{1}\frac{\log\left(1+x^{2}\right)}{1+x^{2}}\,\mathrm{d}x=\pi G-\frac{\pi^{2}}{2}\log2.
	\]
	Putting all things together finally produces 
	\[
	\mathfrak{D}=\pi G-\frac{\pi^{2}}{4}-4\sum_{k=1}^{\infty}\frac{1}{2k+1}\sum_{\ell\geqslant k+1}\frac{\left(-1\right)^{\ell+1}}{\ell}\sum_{m=0}^{\ell-1}\frac{1}{2m+1}.
	\]
	Now since we  know that 
	\[
	\mathfrak{D}=\int_{0}^{\pi/2}\int_{0}^{\pi/2}\dfrac{\log\left(\cos\left(\dfrac{x}{2}\right)\right)-\log\left(\cos\left(\dfrac{z}{2}\right)\right)}{\cos x-\cos z}\,\mathrm{d}x\,\mathrm{d}z=\pi G-\dfrac{7}{4}\,\zeta(3),
	\]
	we finally deduce that
	\[
	\sum_{k=1}^{\infty}\frac{1}{2k+1}\sum_{\ell\geqslant k+1}\frac{\left(-1\right)^{\ell+1}}{\ell}\sum_{m=0}^{\ell-1}\frac{1}{2m+1}=\frac{7}{16}\,\zeta(3)-\dfrac{3}{8}\,\zeta(2),
	\]
	as desired. This completes the proof of Proposition \ref{prop:Prop2}. \QED
	\subsection{Proof of Theorem \ref{thm:generating_function}}
	
	The proof is obtained by substituting the Fourier expansion \ref{eq:Fourier expansion}
	in the intetegral \ref{eq:general integral} and using the identities
	\eqref{eq:by-product1} and \eqref{eq:by-product2}.  \QED
	
	\subsection{Proof of Claim \ref{claim}}
	\label{subsec:Proof-of-Claim}
	
	Replacing $\alpha_{k}\mapsto\alpha^{k}$ in Theorem \ref{thm:generating_function} produces
	\[
	\mathfrak{F}\left(\cos\vartheta\right)=\sum_{k=1}^{\infty}\cos\left(\vartheta\right)\alpha^{k}=\dfrac{\alpha\cos\vartheta-\alpha^{2}}{1-2\alpha\cos\vartheta+\alpha^{2}}.
	\]
	With this choice of the function $\mathfrak{F}$, we have
	\begin{equation}
		\int_{0}^{\pi/2}\int_{0}^{\pi/2}\dfrac{\mathfrak{F}\left(\cos x\right)-\mathfrak{F}\left(\cos z\right)}{\cos x-\cos z}\,\mathrm{d}x\,\mathrm{d}z=\dfrac{4\alpha}{1-\alpha^{2}}\arctan^{2}\left(\dfrac{1+\alpha}{1-\alpha}\right).\label{eq:4alpha}
	\end{equation}
	This can be easily checked as follows: the substitution $\omega=\left(\frac{1+\alpha}{1-\alpha}\right)$
	produces 
	\[
	\mathfrak{F}\left(\cos x\right)=\dfrac{\omega-1}{2}\dfrac{\omega\left(\cos x-1\right)+\left(\cos x+1\right)}{\left(1+\cos x\right)+\omega^{2}\left(1-\cos x\right)}.
	\]
	Therefore, it suffices to compute the following double integral 
	\[
	\int_{0}^{\pi/2}\int_{0}^{\pi/2}\dfrac{1}{\cos x-\cos z}\left\{ \dfrac{\omega\left(\cos x-1\right)+\left(\cos x+1\right)}{\left(1+\cos x\right)+\omega^{2}\left(1-\cos x\right)}-\dfrac{\omega\left(\cos z-1\right)+\left(\cos z+1\right)}{\left(1+\cos z\right)+\omega^{2}\left(1-\cos z\right)}\right\} \mathrm{d}x\,\mathrm{d}z
	\]
	\[
	=\int_{0}^{\pi/2}\int_{0}^{\pi/2}\dfrac{1}{\cos x-\cos z}\left\{ \dfrac{\cos\left(x\right)+1}{(1+\cos x)+\omega^{2}(1-\cos x)}-\dfrac{\cos\left(z\right)+1}{(1+\cos z)+\omega^{2}(1-\cos z)}\right\} \mathrm{d}x\,\mathrm{d}z
	\]
	\[
	+\,\int_{0}^{\pi/2}\int_{0}^{\pi/2}\dfrac{1}{\cos x-\cos z}\left\{ \dfrac{\omega\left(\cos\left(x\right)-1\right)}{(1+\cos x)+\omega^{2}(1-\cos x)}-\dfrac{\omega\left(\cos\left(x\right)-1\right)}{(1+\cos z)+\omega^{2}(1-\cos z)}\right\} \mathrm{d}x\,\mathrm{d}z.
	\]
	Let us first look at the first integral
	\begin{align*}
		&       \int_{0}^{\pi/2}\int_{0}^{\pi/2}\dfrac{1}{\cos x-\cos z}\left\{ \dfrac{\cos\left(x\right)+1}{(1+\cos x)+\omega^{2}(1-\cos x)}-\dfrac{\cos\left(z\right)+1}{(1+\cos z)+\omega^{2}(1-\cos z)}\right\} \mathrm{d}x\,\mathrm{d}z
		\\&=\int_{0}^{\pi/2}\int_{0}^{\pi/2}\dfrac{1}{\cos x-\cos z}\left\{ \dfrac{1}{1+\omega^{2}\left(\frac{1-\cos x}{1+\cos x}\right)}-\dfrac{1}{1+\omega^{2}\left(\frac{1-\cos z}{1+\cos z}\right)}\right\} \mathrm{d}x\,\mathrm{d}z. 
	\end{align*}
	The changes of variable ${\displaystyle \omega_{1}=\tan\left(\dfrac{x}{2}\right)}$
	and $\omega_{2}=\tan\left(\dfrac{z}{2}\right)$ produces
	\[
	\cos\left(x\right)=\left(\dfrac{1-\omega_{1}^{2}}{1+\omega_{1}^{2}}\right),\quad\cos\left(z\right)=\left(\dfrac{1-\omega_{2}^{2}}{1+\omega_{2}^{2}}\right).
	\]
	Therefore, we have
	\begin{align*}
		&\int_{0}^{\pi/2}\int_{0}^{\pi/2}\dfrac{1}{\cos x-\cos z}\left\{ \dfrac{\cos\left(x\right)+1}{(1+\cos x)+\omega^{2}(1-\cos x)}-\dfrac{\cos\left(z\right)+1}{(1+\cos z)+\omega^{2}(1-\cos z)}\right\} \mathrm{d}x\,\mathrm{d}z
		\\&=\int_{0}^{1}\int_{0}^{1}\left(\dfrac{1-\omega_{1}^{2}}{1+\omega_{1}^{2}}-\dfrac{1-\omega_{2}^{2}}{1+\omega_{2}^{2}}\right)^{-1}\left\{ \dfrac{1}{1+\omega^{2}\omega_{1}^{2}}-\dfrac{1}{1+\omega^{2}\omega_{2}^{2}}\right\} \dfrac{4}{\left(1+\omega_{1}^{2}\right)\left(1+\omega_{2}^{2}\right)}\,\mathrm{d}\omega_{1}\mathrm{d}\omega_{2}
		\\&=\int_{0}^{1}\int_{0}^{1}\dfrac{\left(1+\omega_{1}^{2}\right)\left(1+\omega_{2}^{2}\right)}{2\left(\omega_{2}^{2}-\omega_{1}^{2}\right)}\,\dfrac{\omega^{2}\left(\omega_{2}^{2}-\omega_{1}^{2}\right)}{\left(1+\omega^{2}\omega_{1}^{2}\right)\left(1+\omega^{2}\omega_{2}^{2}\right)}\dfrac{4}{\left(1+\omega_{1}^{2}\right)\left(1+\omega_{2}^{2}\right)}\,\mathrm{d}\omega_{1}\mathrm{d}\omega_{2}
		\\&=2\int_{0}^{1}\int_{0}^{1}\dfrac{\omega^{2}}{\left(1+\omega^{2}\omega_{1}^{2}\right)\left(1+\omega^{2}\omega_{2}^{2}\right)}\,\mathrm{d}\omega_{1}\mathrm{d}\omega_{2}=2\arctan^{2}\omega.
	\end{align*}
	Similarly, we find that 
	\begin{align*}
		&\int_{0}^{\pi/2}\int_{0}^{\pi/2}\dfrac{1}{\cos x-\cos z}\left\{ \dfrac{\omega\left(\cos\left(x\right)-1\right)}{(1+\cos x)+\omega^{2}(1-\cos x)}-\dfrac{\omega\left(\cos\left(x\right)-1\right)}{(1+\cos z)+\omega^{2}(1-\cos z)}\right\} \mathrm{d}x\,\mathrm{d}z
		\\&=\omega\int_{0}^{1}\int_{0}^{1}\dfrac{\left(1+\omega_{1}^{2}\right)\left(1+\omega_{2}^{2}\right)}{2\left(\omega_{2}^{2}-\omega_{1}^{2}\right)}\,\left\{ \dfrac{\omega_{1}^{2}}{\omega^{2}\omega_{1}^{2}+1}-\dfrac{\omega_{2}^{2}}{\omega^{2}\omega_{2}^{2}+1}\right\} \dfrac{4}{\left(1+\omega_{1}^{2}\right)\left(1+\omega_{2}^{2}\right)}\,\mathrm{d}\omega_{1}\mathrm{d}\omega_{2}
		\\&     =\omega\int_{0}^{1}\int_{0}^{1}\dfrac{\left(1+\omega_{1}^{2}\right)\left(1+\omega_{2}^{2}\right)}{2\left(\omega_{2}^{2}-\omega_{1}^{2}\right)}\dfrac{\omega_{1}^{2}-\omega_{2}^{2}}{\left(\omega^{2}\omega_{1}^{2}+1\right)\left(\omega^{2}\omega_{2}^{2}+1\right)}\dfrac{4}{\left(1+\omega_{1}^{2}\right)\left(1+\omega_{2}^{2}\right)}\,\mathrm{d}\omega_{1}\mathrm{d}\omega_{2}
		\\&=2\int_{0}^{1}\int_{0}^{1}\dfrac{\omega}{\left(1+\omega^{2}\omega_{1}^{2}\right)\left(1+\omega^{2}\omega_{2}^{2}\right)}\,\mathrm{d}\omega_{1}\mathrm{d}\omega_{2}=\dfrac{2\arctan^{2}\omega}{\omega}.
	\end{align*}
	Putting all things together we recover \eqref{eq:4alpha} and thus obtain the double integral as 
	\begin{align*}
		&\int_{0}^{\pi/2}\int_{0}^{\pi/2}\dfrac{\log\cos\left(\dfrac{x}{2}\right)-\log\cos\left(\dfrac{z}{2}\right)}{\cos x-\cos z}\,\mathrm{d}x\,\mathrm{d}z
		\\&=\sum_{k=1}^{\infty}\dfrac{(-1)^{k-1}}{k}\int_{0}^{\pi/2}\int_{0}^{\pi/2}\dfrac{\cos\left(kx\right)-\cos\left(kz\right)}{\cos x-\cos z}\,\mathrm{d}x\,\mathrm{d}z
		\\&=\sum_{k=1}^{\infty}(-1)^{k-1}\int_{0}^{1}\omega^{k-1}\,\mathrm{d}\omega\int_{0}^{\pi/2}\int_{0}^{\pi/2}\dfrac{\cos\left(kx\right)-\cos\left(kz\right)}{\cos x-\cos z}\,\mathrm{d}x\,\mathrm{d}z
		\\&=\int_{0}^{1}\sum_{k=1}^{\infty}\left(-\omega\right)^{k-1}\mathrm{d}\omega\int_{0}^{\pi/2}\int_{0}^{\pi/2}\dfrac{\cos\left(kx\right)-\cos\left(kz\right)}{\cos x-\cos z}\mathrm{d}x\mathrm{d}z\\&=4\int_{0}^{1}\dfrac{1}{1-\omega^{2}}\arctan^{2}\left(\dfrac{1+\omega}{1-\omega}\right)\mathrm{d}\omega
		\\&=4\int_{0}^{1}\dfrac{1}{1-\omega^{2}}\left(\dfrac{\pi}{4}-\arctan\omega\right)^{2}\mathrm{d}\omega
		\\&=4\int_{0}^{\pi/4}\dfrac{1}{\cos\left(2\vartheta\right)}\left(\dfrac{\pi}{4}-\vartheta\right)^{2}\mathrm{d}\vartheta\\&=4\int_{0}^{\pi/4}\dfrac{\vartheta^{2}}{\sin\left(2\vartheta\right)}\,\mathrm{d}\vartheta\\&=\pi G-\dfrac{7}{4}\,\zeta(3),
	\end{align*}
	which is the desired result. \QED
	
	%
	%
	%
	
	
	\subsection{Proof of Theorem \ref{thm:representations}}
	In a recent preprint \cite{Xu} by C. Xu, we find the identity 
	\[
	\int_{x<y<z<1}\frac{\mathrm{d}y\,\mathrm{d}z}{\left(1+y^{2}\right)\left(1+z^{2}\right)}=\frac{1}{2}\left(\frac{\pi}{4}-\arctan x\right)^{2},
	\]
	from which we deduce the three-dimensional representation 
	\begin{align*}
		&\frac{1}{8}\int_{0}^{\pi/2}\int_{0}^{\pi/2}\dfrac{\log\left(\cos\left(\dfrac{x}{2}\right)\right)-\log\left(\cos\left(\dfrac{z}{2}\right)\right)}{\cos x-\cos z}\,\mathrm{d}x\,\mathrm{d}z  \\&=\frac{1}{2}\int_{0}^{1}\frac{1}{1-x^{2}}\arctan^{2}\left(\dfrac{1-x}{1+x}\right)\mathrm{d}x
		\\&=\frac{1}{2}\int_{0}^{1}\frac{\left(\pi/4-\arctan x\right)^{2}}{1-x^{2}}\,\mathrm{d}x
		\\&=\int_{0<x<y<z<1}\frac{\mathrm{d}x\,\mathrm{d}y\,\mathrm{d}z}{\left(1-x^{2}\right)\left(1+y^{2}\right)\left(1+z^{2}\right)}.
	\end{align*}
	This identity produces a MZV representation of our integral as follows:
	since
	\[
	\int_{0}^{y}\frac{1}{1-x^{2}}\,\mathrm{d}x=\sum_{k=0}^{\infty}\frac{y^{2k+1}}{2k+1},
	\]
	\begin{align*}
		\int_{0}^{z}\sum_{k=0}^{\infty}\frac{1}{1+y^{2}}\frac{y^{2k+1}}{2k+1}\,\mathrm{d}y&=\sum_{k=0}^{\infty}\sum_{l=0}^{\infty}\left(-1\right)^{l}\int_{0}^{z}y^{2l}\,\frac{y^{2k+1}}{2k+1}\,\mathrm{d}y\\&=\sum_{k=0}^{\infty}\sum_{l=0}^{\infty}\frac{\left(-1\right)^{l}}{2k+1}\frac{z^{2k+2l+2}}{2k+2l+2},
	\end{align*}
	and finally we have
	\begin{align*}
		&\int_{0}^{1}\frac{\mathrm{d}z}{1+z^{2}}\sum_{k,l,m\geqslant0}
		\frac{\left(-1\right)^{l}}{2k+1}\frac{z^{2k+2l+2}}{2k+2l+2} \\&=\sum_{k,l,m\geqslant0}\frac{\left(-1\right)^{l+m}}{\left(2k+1\right)\left(2k+2l+2\right)}\int_{0}^{1}z^{2k+2l+2m+2}\,\mathrm{d}z\\ &=\sum_{k,l,m\geqslant0}\frac{\left(-1\right)^{l+m}}{\left(2k+1\right)\left(2k+2l+2\right)\left(2k+2l+2m+3\right)}.\end{align*}
	Therefore, we deduce that
	\begin{align*}&\int_{0<x<y<z<1}\frac{\mathrm{d}x\,\mathrm{d}y\,\mathrm{d}z}{\left(1-x^{2}\right)\left(1+y^{2}\right)\left(1+z^{2}\right)}\\&=\sum_{k,l,m\geqslant0}\frac{\left(-1\right)^{l+m}}{\left(2k+1\right)\left(2k+2l+2\right)\left(2k+2l+2m+3\right)}\\& =\sum_{k,l,m\geqslant0}\frac{\left(-1\right)^{l+m}}{\left(2k+1\right)\left(2k+1+2l+1\right)\left(2k+1+2l+1+2m+1\right)}\\&=\sum_{\underset{\text{all odd}}{0<k<l<m}}\frac{\left(-1\right)^{\frac{m+k}{2}-1}}{klm},
	\end{align*}
	which is the desired result. \QED
	\subsection{\label{subsec:Proof-of-Proposition 8}Proof of Proposition \ref{catalan_triple}}
	First, let us break down the triple integral as follows
	\begin{align*}
		&8\int_{0}^{1}\int_{0}^{\vartheta_{3}}\int_{0}^{\vartheta_{2}}\dfrac{\mathrm{d}\vartheta_{1}\,\mathrm{d}\vartheta_{2}\,\mathrm{d}\vartheta_{3}}{\left(1-\vartheta_{1}^{2}\right)\left(1+\vartheta_{2}^{2}\right)\left(1+\vartheta_{3}^{2}\right)}\\&=-4\int_{0}^{1}\dfrac{1}{1+\vartheta_{3}^{2}}\int_{0}^{\vartheta_{3}}\dfrac{1}{1+\vartheta_{2}^{2}}\log\left(\frac{1-\vartheta_{2}}{1+\vartheta_{2}}\right)\mathrm{d}\vartheta_{2}\,\mathrm{d}\vartheta_{3}
		\\&=-\,4\int_{0}^{1}\dfrac{1}{1+\vartheta_{2}^{2}}\log\left(\frac{1-\vartheta_{2}}{1+\vartheta_{2}}\right)\int_{\vartheta_{2}}^{1}\dfrac{1}{1+\vartheta_{3}^{2}}\,\mathrm{d}\vartheta_{3}\,\mathrm{d}\vartheta_{2}\\&\,\,\,\,\,\,\,-4\int_{0}^{1}\dfrac{1}{1+\vartheta_{2}^{2}}\left(\dfrac{\pi}{4}-\arctan\vartheta_{2}\right)\log\left(\frac{1-\vartheta_{2}}{1+\vartheta_{2}}\right)\mathrm{d}\vartheta_{2}.
	\end{align*}
	Since 
	\[
	\int_{0}^{1}\dfrac{1}{1+\vartheta_{2}^{2}}\log\left(\frac{1-\vartheta_{2}}{1+\vartheta_{2}}\right)\int_{\vartheta_{2}}^{1}\dfrac{1}{1+\vartheta_{3}^{2}}\,\mathrm{d}\vartheta_{3}\,\mathrm{d}\vartheta_{2}=0,
	\]
	and $\frac{\pi}{4}-\arctan\vartheta=\arctan\left(\frac{1-\vartheta}{1+\vartheta}\right)$, we deduce that 
	\[
	8\int_{0}^{1}\int_{0}^{\vartheta_{3}}\int_{0}^{\vartheta_{2}}\dfrac{\mathrm{d}\vartheta_{1}\,\mathrm{d}\vartheta_{2}\,\mathrm{d}\vartheta_{3}}{\left(1-\vartheta_{1}^{2}\right)\left(1+\vartheta_{2}^{2}\right)\left(1+\vartheta_{3}^{2}\right)}=-4\int_{0}^{1}\dfrac{1}{1+\vartheta_{2}^{2}}\arctan\left(\frac{1-\vartheta_{2}}{1+\vartheta_{2}}\right)\log\left(\frac{1-\vartheta_{2}}{1+\vartheta_{2}}\right)\mathrm{d}\vartheta_{2}.
	\]
	Substituting $\left(\frac{1-\vartheta_{2}}{1+\vartheta_{2}}\right)\mapsto\vartheta_{2}$ yields
	\begin{align*}
		&       -4\int_{0}^{1}\dfrac{1}{1+\vartheta_{2}^{2}}\arctan\left(\frac{1-\vartheta_{2}}{1+\vartheta_{2}}\right)\log\left(\frac{1-\vartheta_{2}}{1+\vartheta_{2}}\right)\mathrm{d}\vartheta_{2}\\&=-4\int_{0}^{1}\dfrac{\arctan\vartheta_{2}\log\vartheta_{2}}{1+\vartheta_{2}^{2}}\,\mathrm{d}\vartheta_{2}
		\\&=-2\int_{0}^{\infty}\dfrac{\arctan\vartheta_{2}\log\vartheta_{2}}{1+\vartheta_{2}^{2}}\,\mathrm{d}\vartheta_{2}-\pi\int_{0}^{1}\dfrac{\log\vartheta_{2}}{1+\vartheta_{2}^{2}}\,\mathrm{d}\vartheta_{2}\\&=-2\int_{0}^{\pi/2}\vartheta\log\left(\tan\vartheta\right)\mathrm{d}\vartheta+\pi G.
	\end{align*}
	Next, we use the Fourier series expansions 
	\[
	\log\left(\sin\vartheta\right)=\sum_{n=1}^{\infty}(-1)^{n}\,\frac{\cos\left(2n\vartheta\right)}{n},
	\]
	and identity \eqref{eq:Fourier logcos} for $\log\left(\cos\vartheta\right),$
	to finally obtain
	\begin{align*}
		\int_{0}^{\pi/2}\vartheta\log\left(\tan\vartheta\right)\,\mathrm{d}\vartheta&=\int_{0}^{\pi/2}\vartheta\left(\sum_{n=1}^{\infty}(-1)^{n}\,\frac{\cos\left(2n\vartheta\right)}{n}-\sum_{n=1}^{\infty}\frac{\cos\left(2n\vartheta\right)}{n}\right)\mathrm{d}\vartheta
		\\&=\sum_{n=1}^{\infty}\frac{(-1)^{n}}{n}\left[\frac{\pi\sin\left(n\pi\right)}{4n}+\frac{\cos\left(n\pi\right)}{4n^{2}}-\frac{1}{4n^{2}}\right]\\&-\sum_{n=1}^{\infty}\frac{1}{n}\left[\frac{\pi\sin\left(n\pi\right)}{4n}+\frac{\cos\left(n\pi\right)}{4n^{2}}-\frac{1}{4n^{2}}\right]=\dfrac{7}{8}\,\zeta(3).
	\end{align*}
	Putting all things together gives us the desired result.\QED
	\subsection{Proof of proposition \ref{lol_int}}
	Start with the identity 
	\[
	\int_{-1}^{1}\frac{1}{1+xz}\frac{\mathrm{d}x}{\sqrt{1-x^{2}}}=\frac{\pi}{\sqrt{1-z^{2}}}.
	\]
	Integrating with respect to $z$ on both sides
	produces
	\[
	\int\frac{\pi^{2}}{1-z^{2}}\,\mathrm{d}z=\frac{\pi^{2}}{2}\log\left(\frac{1+z}{1-z}\right),
	\]
	on the right-hand side and \begin{align*}
		&\int\mathrm{d}z\int_{-1}^{1}\int_{-1}^{1}\frac{1}{x_{1}-x_{2}}\left[\frac{x_{1}}{1+x_{1}z}-\frac{x_{2}}{1+x_{2}z}\right]\frac{\mathrm{d}x_{1}}{\sqrt{1-x_{1}^{2}}}\frac{\mathrm{d}x_{2}}{\sqrt{1-x_{2}^{2}}}\\&=\int_{-1}^{1}\int_{-1}^{1}\frac{\log\left(1+zx_{1}\right)-\log\left(1+zx_{2}\right)}{x_{1}-x_{2}}\frac{\mathrm{d}x_{1}}{\sqrt{1-x_{1}^{2}}}\frac{\mathrm{d}x_{2}}{\sqrt{1-x_{2}^{2}}}\\&=\int_{0}^{\pi}\int_{0}^{\pi}\frac{\log\left(1+z\cos\left(\vartheta_{1}\right)\right)-\log\left(1+z\cos\left(\vartheta_{2}\right)\right)}{\cos\vartheta_{1}-\cos\vartheta_{2}}\,\mathrm{d}\vartheta_{1}\,\mathrm{d}\vartheta_{2},
	\end{align*}
	on the left-hand side, so that we have 
	\[
	\int_{0}^{\pi}\int_{0}^{\pi}\frac{\log\left(1+z\cos\vartheta_{1}\right)-\log\left(1+z\cos\vartheta_{2}\right)}{\cos\vartheta_{1}-\cos\vartheta_{2}}\,\mathrm{d}\vartheta_{1}\,\mathrm{d}\vartheta_{2}=\frac{\pi^{2}}{2}\log\left(\frac{1+z}{1-z}\right)+C.
	\]
	The constant is evaluated to $C=0$ by choosing $z=0$ on
	both sides. Thus, we have 
	\[
	\int_{0}^{\pi}\int_{0}^{\pi}\frac{\log\left(1+z\cos\vartheta_{1}\right)-\log\left(1+z\cos\vartheta_{2}\right)}{\cos\vartheta_{1}-\cos\vartheta_{2}}\,\mathrm{d}\vartheta_{1}\,\mathrm{d}\vartheta_{2}=\frac{\pi^{2}}{2}\log\left(\frac{1+z}{1-z}\right),
	\]
	as desired. \QED
	\subsection{Proof of proposition \ref{polylog_int}}
	To start with, notice that for $k\geqslant2$, we have 
	\begin{align}
		\int_{0}^{\pi}\int_{0}^{\pi}\dfrac{\cos^{k}\vartheta_{1}-\cos^{k}\vartheta_{2}}{\cos\vartheta_{1}-\cos\vartheta_{2}}\,\mathrm{d}\vartheta_{1}\,\mathrm{d}\vartheta_{2}&=\sum_{n=0}^{k-1}\int_{0}^{\pi}\int_{0}^{\pi}\cos^{n}\vartheta_{1}\cos^{n-k-1}\vartheta_{2}\,\mathrm{d}\vartheta_{1}\,\mathrm{d}\vartheta_{2}\nonumber
		\\& =\sum_{n=0}^{k-1}\int_{0}^{\pi}\cos^{n}\vartheta_{1}\mathrm{d}\vartheta_{1}\int_{0}^{\pi}\cos^{n-k-1}\vartheta_{2}\,\mathrm{d}\vartheta_{2}\nonumber\\&=\begin{cases}
			\pi^{2}&k\,\,\text{is odd}.\\
			0&k\,\,\text{is even}.
		\end{cases}\label{cosk_result}
	\end{align}
	Next, we deduce that 
	\begin{align*}
		&\int_{0}^{\pi}\int_{0}^{\pi}\dfrac{\mathrm{Li}_{n}\left(z\cos\vartheta_{1}\right)-\mathrm{Li}_{n}\left(z\cos\vartheta_{2}\right)}{\cos\vartheta_{1}-\cos\vartheta_{2}}\mathrm{d}\vartheta_{1}\mathrm{d}\vartheta_{2}=\sum_{k=1}^{\infty}\dfrac{z^{k}}{k^{n}}\int_{0}^{\pi}\int_{0}^{\pi}\dfrac{\cos^{k}\vartheta_{1}-\cos^{k}\vartheta_{2}}{\cos\vartheta_{1}-\cos\vartheta_{2}}\mathrm{d}\vartheta_{1}\mathrm{d}\vartheta_{2}.
	\end{align*}
	Replacing $k\mapsto2k-1$ in the above equation and using equation
	\ref{cosk_result} yields 
	\begin{align*}
		&\sum_{k=1}^{\infty}\dfrac{z^{2k-1}}{(2k-1)^{n}}\int_{0}^{\pi}\int_{0}^{\pi}\dfrac{\cos^{2k-1}\vartheta_{1}-\cos^{2k-1}\vartheta_{2}}{\cos\vartheta_{1}-\cos\vartheta_{2}}\,\mathrm{d}\vartheta_{1}\,\mathrm{d}\vartheta_{2}=\pi^{2}\sum_{k=1}^{\infty}\dfrac{z^{2k-1}}{(2k-1)^{n}}
		\\& =\pi^{2}\left(\sum_{k=1}^{\infty}\dfrac{z^{k}}{k^{n}}-\dfrac{1}{2^{n}}\sum_{k=1}^{\infty}\dfrac{z^{2k}}{k^{n}}\right)=\pi^{2}\left(\Li_{n}\left(z\right)-\dfrac{1}{2^{n}}\Li_{n}\left(z^{2}\right)\right)=\dfrac{\pi^{2}}{2}\left[\mathrm{Li}_{n}\left(z\right)-\mathrm{Li}_{n}\left(-z\right)\right].
	\end{align*}
	Substituting $z=1$ in the previous result produces
	\begin{align*}
		\int_{0}^{\pi}\int_{0}^{\pi}\dfrac{\mathrm{Li}_{n}\left(\cos\vartheta_{1}\right)-\mathrm{Li}_{n}\left(\cos\vartheta_{2}\right)}{\cos\vartheta_{1}-\cos\vartheta_{2}}\mathrm{d}\vartheta_{1}\mathrm{d}\vartheta_{2}&=\pi^{2}\left(\sum_{k=1}^{\infty}\dfrac{1}{k^{n}}-\dfrac{1}{2^{n}}\sum_{k=1}^{\infty}\dfrac{1}{k^{n}}\right)\\&=\pi^{2}\left(1-2^{-n}\right)\zeta(n).
	\end{align*}
	Now substituting $z=\phi$ in Proposition \ref{polylog_int} and using 
	\[
	\mathrm{Li}_{2}\left(\phi\right)=\dfrac{\pi^{2}}{10}-\log^{2}\left(\dfrac{1+\sqrt{5}}{2}\right),\quad
	\mathrm{Li}_{2}\left(-\phi\right)=\dfrac{1}{2}\log^{2}\left(\dfrac{1+\sqrt{5}}{2}\right)-\dfrac{\pi^{2}}{15},
	\]
	we deduce that 
	\begin{align*}
		&\int_{0}^{\pi}\int_{0}^{\pi}\dfrac{\mathrm{Li}_{2}\left(\phi\cos\vartheta_{1}\right)-\mathrm{Li}_{2}\left(\phi\cos\vartheta_{2}\right)}{\cos\vartheta_{1}-\cos\vartheta_{2}}\,\mathrm{d}\vartheta_{1}\,\mathrm{d}\vartheta_{2}=\dfrac{\pi^{2}}{2}\left[\mathrm{Li}_{2}\left(\phi\right)-\mathrm{Li}_{2}\left(-\phi\right)\right]
		\\&=\dfrac{\pi^{2}}{2}\left[\left(\dfrac{\pi^{2}}{10}-\log^{2}\left(\dfrac{1+\sqrt{5}}{2}\right)\right)-\left(\dfrac{1}{2}\log^{2}\left(\dfrac{1+\sqrt{5}}{2}\right)-\dfrac{\pi^{2}}{15}\right)\right]=\dfrac{\pi^{4}}{12}-\dfrac{3}{4}\log^{2}\left(\dfrac{1+\sqrt{5}}{2}\right),
	\end{align*}
	which is the desired result. \QED
	\subsection{Proof of Theorem \ref{analytic_thm}}
	The integral is expanded as 
	\[
	\int_{0}^{\pi}\int_{0}^{\pi}\frac{\mathfrak{F}\left(\cos\vartheta_{1}\right)-\mathfrak{F}\left(\cos\vartheta_{2}\right)}{\cos\vartheta_{1}-\cos\vartheta_{2}}\,\mathrm{d}\vartheta_{1}\,\mathrm{d}\vartheta_{2}=\sum_{n=0}^{\infty}a_{n}z^{n}\int_{0}^{\pi}\int_{0}^{\pi}\frac{\cos^{n}\vartheta_{1}-\cos^{n}\vartheta_{2}}{\cos\vartheta_{1}-\cos\vartheta_{2}}\,\mathrm{d}\vartheta_{1}\,\mathrm{d}\vartheta_{2}.
	\]
	From equation \eqref{cosk_result}, we have
	\[
	\sum_{k=0}^{n-1}\int_{0}^{\pi}\cos^{k}\vartheta_{1}\mathrm{d}\vartheta_{1}\int_{0}^{\pi}\cos^{k-n-1}\vartheta_{2}\,\mathrm{d}\vartheta_{2}=\begin{cases}
		\pi^{2} & n \,\,\text{is even}.\\
		0 & n \,\,\text{is odd}.
	\end{cases}
	\]
	Indeed, from the generating function 
	\begin{align*}
		\sum_{k=0}^{\infty}z^{k}\int_{0}^{\pi}\cos^{k}\vartheta_{1}\mathrm{d}\vartheta_{1}&=\sum_{k=0}^{\infty}\Gamma\left(k+\dfrac{1}{2}\right)\Gamma\left(\dfrac{1}{2}\right)\frac{z^{2k}}{\Gamma\left(k+1\right)}\\&=\pi\sum_{k=0}^{\infty}\left(\frac{1}{2}\right)_{k}\frac{z^{2k}}{k!}=\frac{\pi}{\sqrt{1-z^{2}}},   
	\end{align*}
	it follows that the sum $$\sum_{k=0}^{n-1}\int_{0}^{\pi}\cos^{k}\vartheta_{1}\mathrm{d}\vartheta_{1}\int_{0}^{\pi}\cos^{k-n-1}\vartheta_{2}\,\mathrm{d}\vartheta_{2},$$
	is the coefficient of $z^{n-1}$ in the generating function $\pi^2/\left(1-z^2\right)$.
	
	This result includes the special case $n=0$ as well, so that we have
	\begin{align*}
		\int_{0}^{\pi}\int_{0}^{\pi}\frac{\mathfrak{F}\left(\cos\vartheta_{1}\right)-\mathfrak{F}\left(\cos\vartheta_{2}\right)}{\cos\vartheta_{1}-\cos\vartheta_{2}}\,\mathrm{d}\vartheta_{1}\,\mathrm{d}\vartheta_{2} =\pi^{2}\sum_{n=0}^{\infty}a_{2n+1}z^{2n+1}.
	\end{align*}
	The other case 
	\[
	\int_{0}^{\pi}\int_{0}^{\pi}\frac{\cos\vartheta_{1}\mathfrak{F}\left(\cos\vartheta_{1}\right)-\cos\vartheta_{2}\mathfrak{F}\left(\cos\vartheta_{2}\right)}{\cos\vartheta_{1}-\cos\vartheta_{2}}\,\mathrm{d}\vartheta_{1}\,\mathrm{d}\vartheta_{2},
	\]
	is obtained similarly by expanding 
	\begin{align*}
		&\int_{0}^{\pi}\int_{0}^{\pi}\frac{\cos\vartheta_{1}\mathfrak{F}\left(\cos\vartheta_{1}\right)-\cos\vartheta_{2}\mathfrak{F}\left(\cos\vartheta_{2}\right)}{\cos\vartheta_{1}-\cos\vartheta_{2}}\,\mathrm{d}\vartheta_{1}\,\mathrm{d}\vartheta_{2}\\&=\sum_{n=0}^{\infty}a_{n}z^{n}\int_{0}^{\pi}\int_{0}^{\pi}\frac{\cos^{n+1}\vartheta_{1}-\cos^{n+1}\vartheta_{2}}{\cos\vartheta_{1}-\cos\vartheta_{2}}\,\mathrm{d}\vartheta_{1}\,\mathrm{d}\vartheta_{2},
	\end{align*}
	and remarking that, similarly, 
	\[
	\int_{0}^{\pi}\int_{0}^{\pi}\frac{\cos^{n+1}\vartheta_{1}-\cos^{n+1}\vartheta_{2}}{\cos\vartheta_{1}-\cos\vartheta_{2}}\,\mathrm{d}\vartheta_{1}\,\mathrm{d}\vartheta_{2}=\begin{cases}
		\pi^{2} &  n \,\,\text{is odd}.\\
		0 & n \,\,\text{is even}.
	\end{cases}
	\]
	This completes the proof. \QED
	\subsection{Proof of Proposition \ref{one_parameter}}
	\label{one_parameter_proof}
	Start from the identity 
	\[
	\int_{0}^{1}\frac{1}{1+xz}\frac{\mathrm{d}x}{\sqrt{1-x^{2}}}=\frac{\arccos z}{\sqrt{1-z^{2}}},
	\]
	so that we have 
	\[
	\int\mathrm{d}z\int_{0}^{1}\int_{0}^{1}\frac{1}{1+x_{1}z}\frac{1}{1+x_{2}z}\frac{\mathrm{d}x_{1}}{\sqrt{1-x_{1}^{2}}}\frac{\mathrm{d}x_{2}}{\sqrt{1-x_{2}^{2}}}=\int\left(\frac{\arccos z}{\sqrt{1-z^{2}}}\right)^{2}\mathrm{d}z.
	\]
	Next, we notice that 
	\begin{align*}&
		\int\mathrm{d}z\int_{0}^{1}\int_{0}^{1}\frac{1}{1+x_{1}z}\frac{1}{1+x_{2}z}\frac{\mathrm{d}x_{1}}{\sqrt{1-x_{1}^{2}}}\frac{\mathrm{d}x_{2}}{\sqrt{1-x_{2}^{2}}}
		\\&=\int_{0}^{1}\int_{0}^{1}\frac{\log\left(1+x_{1}z\right)-\log\left(1+x_{2}z\right)}{x_{1}-x_{2}}\frac{\mathrm{d}x_{1}}{\sqrt{1-x_{1}^{2}}}\frac{\mathrm{d}x_{2}}{\sqrt{1-x_{2}^{2}}}
		\\&
		=\int_{0}^{\pi/2}\int_{0}^{\pi/2}\frac{\log\left(1+z\cos\vartheta_{1}\right)-\log\left(1+z\cos\vartheta_{2}\right)}{\cos\vartheta_{1}-\cos\vartheta_{2}}\,\mathrm{d}\vartheta_{1}\,\mathrm{d}\vartheta_{2},
	\end{align*}
	where the last equality is deduced from the changes of variable $x_{1}=\cos\vartheta_{1}$
	and $x_{2}=\cos\vartheta_{2}$. 
	
	With the change of variable $z=\cos\omega,$
	the right-hand side reads 
	\begin{align*}
		&\int\left(\frac{\arccos\left(z\right)}{\sqrt{1-z^{2}}}\right)^{2}\mathrm{d}z=\int\frac{\arccos^{2}z}{1-z^{2}}\,\mathrm{d}z=\int\dfrac{\omega^{2}}{\sin\omega}\,\mathrm{d}\omega
		\\& =2\imath\omega\left[\Li_{2}\left(e^{\imath\omega}\right)-\Li_{2}\left(-e^{\imath\omega}\right)\right]+2\left[\Li_{3}\left(-e^{\imath\omega}\right)-\Li_{3}\left(-e^{\imath\omega}\right)\right]-\omega^{2}\left(\dfrac{1-e^{\imath\omega}}{1+e^{\imath\omega}}\right).
	\end{align*}
	Putting all things together produces 
	\begin{align*}
		&\int_{0}^{\pi/2}\int_{0}^{\pi/2}\frac{\log\left(1+z\cos\vartheta_{1}\right)-\log\left(1+z\cos\vartheta_{2}\right)}{\cos\vartheta_{1}-\cos\vartheta_{2}}\,\mathrm{d}\vartheta_{1}\,\mathrm{d}\vartheta_{2}
		\\& =2\imath\arccos z\left[\Li_{2}\left(e^{\imath\arccos z}\right)-\Li_{2}\left(-e^{\imath\arccos z}\right)\right]+2\left[\Li_{3}\left(-e^{\imath\arccos z}\right)-\Li_{3}\left(e^{\imath\arccos z}\right)\right]
		\\& -\arccos^{2}z\log\left(\dfrac{1-e^{\imath\arccos z}}{1+e^{\imath\arccos z}}\right)+K,
	\end{align*}
	so that, by replacing $z$ by $\cos z$, we have
	\begin{align*}
		& \int_{0}^{\pi/2}\int_{0}^{\pi/2}\frac{\log\left(1+\cos z\cos\vartheta_{1}\right)-\log\left(1+\cos z\cos\vartheta_{2}\right)}{\cos\vartheta_{1}-\cos\vartheta_{2}}\,\mathrm{d}\vartheta_{1}\,\mathrm{d}\vartheta_{2}
		\\& =2\imath z\left[\Li_{2}\left(e^{\imath z}\right)-\Li_{2}\left(-e^{\imath z}\right)\right]+2\left[\Li_{3}\left(-e^{\imath z}\right)-\Li_{3}\left(e^{\imath z}\right)\right] -z^{2}\log\left(\dfrac{1-e^{\imath z}}{1+e^{\imath z}}\right)+K.
	\end{align*}
	The integration constant $K$ is evaluated to by choosing $z=\pi/2$
	on both sides to produce
	\[
	0=K-2\pi G.
	\]
	This is because at $z=\pi/2$, we have 
	\[
	\mathrm{Li}_{2}\left(\imath\right)=\imath G-\dfrac{\pi^{2}}{48},\quad\mathrm{Li}_{2}\left(-\imath\right)=-\imath G-\dfrac{\pi^{2}}{48},
	\]
	so that 
	\[
	2\imath z\left[\mathrm{Li}_{2}\left(e^{\imath z}\right)-\mathrm{Li}_{2}\left(-e^{\imath z}\right)\right]=-2\pi G.
	\]
	Moreover, we also find that 
	\[
	\mathrm{Li}_{3}\left(\imath\right)-\mathrm{Li}_{3}\left(-\imath\right)=\dfrac{\imath\pi^{3}}{16},
	\]
	whereas, at $z=\pi/2$, the last term evaluates as 
	\[
	z^{2}\log\left(\frac{1-e^{\imath z}}{1+e^{\imath z}}\right)=\left(\frac{\pi}{2}\right)^{2}\log\left(\frac{1-\imath}{1+\imath}\right)=\left(\frac{\pi}{2}\right)^{2}\left(-\frac{\imath\pi}{2}\right)=-\frac{\imath\pi^{3}}{8},
	\]
	so that finally we have, at $z=\pi/2$,
	\[
	2\left[\mathrm{Li}_{3}\left(-e^{\imath z}\right)-\mathrm{Li}_{3}\left(e^{\imath z}\right)\right]-z^{2}\log\left(\frac{1-e^{\imath z}}{1+e^{\imath z}}\right)=0.
	\]
	For $z=0$, the only term that remains
	is
	\[
	2\left(\Li_{3}\left(-1\right)-\Li_{3}\left(1\right)\right)+2\pi G=2\pi G-\dfrac{7}{2}\,\zeta(3).
	\]
	This completes the proof. \QED
	\subsection{Proof of Theorem \ref{lil_general}}
	First, we write
	\begin{align*}
		&\int_{0}^{\pi/2}\int_{0}^{\pi/2}\frac{\mathfrak{F}\left(z\cos\vartheta_{1}\right)-\mathfrak{F}\left(z\cos\vartheta_{2}\right)}{\cos\vartheta_{1}-\cos\vartheta_{2}}\,\mathrm{d}\vartheta_{1}\,\mathrm{d}\vartheta_{2}\\&=\sum_{n=1}^{\infty}\mathfrak{R}_{n}z^{n}\int_{0}^{\pi/2}\int_{0}^{\pi/2} \frac{\cos^{n}\vartheta_{1}-\cos^{n}\vartheta_{2}}{\cos\vartheta_{1}-\cos\vartheta_{2}}\,\mathrm{d}\vartheta_{1}\,\mathrm{d}\vartheta_{2}.
	\end{align*}
	Notice that the coefficient of $\mathfrak{R}_{0}$ in the right-hand side is
	$0,$ hence the sum starts at $n=1.$

	Next we write, for $n\geqslant1,$
	\[
	\frac{\cos^{n}\vartheta_{1}-\cos^{n}\vartheta_{2}}{\cos\vartheta_{1}-\cos\vartheta_{2}}=\sum_{k=0}^{n-1}\cos^{k}\vartheta_{1}\cos^{n-k-1}\vartheta_{2},
	\]
	so that we have
	\[
	\int_{0}^{\pi/2}\int_{0}^{\pi/2}\frac{\cos^{n}\vartheta_{1}-\cos^{n}\vartheta_{2}}{\cos\vartheta_{1}-\cos\vartheta_{2}}\,\mathrm{d}\vartheta_{1}\,\mathrm{d}\vartheta_{2}=\sum_{k=0}^{n-1}\int_{0}^{\pi/2}\cos^{k}\vartheta_{1}\mathrm{d}\vartheta_{1}\int_{0}^{\pi/2}\cos^{n-k-1}\vartheta_{2}\,\mathrm{d}\vartheta_{2}.
	\]
	With
	\[
	\int_{0}^{\pi/2}\cos^{k}\vartheta\, \mathrm{d}\vartheta=\frac{\sqrt{\pi}}{2}\,\Gamma\left(\frac{k+1}{2}\right)\Gamma\left(\frac{k}{2}+1\right)^{-1},
	\]
	we have
	\begin{align*}
		&\int_{0}^{\pi/2}\int_{0}^{\pi/2}\frac{\cos^{n}\vartheta_{1}-\cos^{n}\vartheta_{2}}{\cos\vartheta_{1}-\cos\vartheta_{2}}\,\mathrm{d}\vartheta_{1}\,\mathrm{d}\vartheta_{2}\\&=\frac{\pi}{4}\sum_{k=0}^{n-1}\Gamma\left(\frac{k+1}{2}\right)\Gamma\left(\frac{k}{2}+1\right)^{-1}\Gamma\left(\frac{n-k}{2}\right)\Gamma\left(\frac{n-k-1}{2}+1\right)^{-1}.
	\end{align*}
	But this sum is the coefficient $\tau_{n-1}$ of $\omega^{n-1}$ in the generating
	function
	\[
	\left(\sum_{k=0}^{\infty}\Gamma\left(\frac{k+1}{2}\right)\Gamma\left(\frac{k}{2}+1\right)^{-1}\omega^{k}\right)^{2}=\left(\frac{\arccos\sqrt{1-\omega^{2}}}{\sqrt{1-\omega^{2}}}\right)^{2}.
	\]
	Therefore, we deduce that
	\[
	\int_{0}^{\pi/2}\int_{0}^{\pi/2}\frac{\mathfrak{F}\left(z\cos\vartheta_{1}\right)-\mathfrak{F}\left(z\cos\vartheta_{2}\right)}{\cos\vartheta_{1}-\cos\vartheta_{2}}\,\mathrm{d}\vartheta_{1}\,\mathrm{d}\vartheta_{2}=\sum_{n=1}^{\infty}\mathfrak{R}_{n}\mathfrak{H}_{n}z^{n},
	\]
	with $\mathfrak{H}_{n}=\tau_{n-1}$ and
	\[
	\left(\frac{\arccos\sqrt{1-\omega^{2}}}{\sqrt{1-\omega^{2}}}\right)^{2}=\sum_{n=2}^{\infty}\tau_{n}\,\omega^{n},
	\]
	so that we have
	\begin{align*}
		\sum_{n=1}^{\infty}\tau_{n-1}\,\omega^{n}=\sum_{n=0}^{\infty}\tau_{n}\,\omega^{n+1}=\omega\left(\frac{\arccos\sqrt{1-\omega^{2}}}{\sqrt{1-\omega^{2}}}\right)^{2}.
	\end{align*}
	Next, we use the integral formula for the Hadamard product of $$\mathfrak{F}\left(z\right)=\sum_{n=0}^{\infty}\mathfrak{R}_{n}z^{n},\quad G\left(z\right)=\sum_{n=0}^{\infty}\mathfrak{H}_{n}z^{n},$$
	to get
	\[
	\left(\mathfrak{F}\otimes G\right)\left(z\right)=\sum_{n=1}^{\infty}\mathfrak{R}_{n}\mathfrak{H}_{n}z^{n}=\frac{1}{2\pi}\int_{0}^{2\pi}\mathfrak{F}\left(e^{\imath\vartheta}\sqrt{z}\right)G\left(e^{-\imath\vartheta}\sqrt{z}\right)\mathrm{d}\vartheta,
	\]
	and we deduce that
	\begin{align*}
		&\int_{0}^{\pi/2}\int_{0}^{\pi/2}\frac{\mathfrak{F}\left(z\cos\vartheta_{1}\right)-\mathfrak{F}\left(z\cos\vartheta_{2}\right)}{\cos\vartheta_{1}-\cos\vartheta_{2}}\,\mathrm{d}\vartheta_{1}\,\mathrm{d}\vartheta_{2}
		\\&=\frac{1}{2\pi}\int_{0}^{2\pi}\mathfrak{F}\left(e^{\imath\vartheta}\sqrt{z}\right)\frac{e^{-\imath\vartheta}\sqrt{z}}{1-ze^{2\imath\vartheta}}\arccos^{2}\left(\sqrt{1-ze^{2\imath\vartheta}}\right)\mathrm{d}\vartheta.\end{align*}
	As a special case $z=1$, we obtain
	\[
	\int_{0}^{\pi/2}\int_{0}^{\pi/2}\frac{\mathfrak{F}\left(\cos\vartheta_{1}\right)-\mathfrak{F}\left(\cos\vartheta_{2}\right)}{\cos\vartheta_{1}-\cos\vartheta_{2}}\,\mathrm{d}\vartheta_{1}\,\mathrm{d}\vartheta_{2}=\frac{1}{4\pi \imath}\int_{0}^{2\pi}\frac{\mathfrak{F}\left(e^{\imath\vartheta}\right)}{\sin\left(\vartheta\right)}\arccos^{2}\left(\sqrt{1-e^{-2\imath\vartheta}}\right)\mathrm{d}\vartheta.
	\]
	as desired. \QED
	\subsection{Proof of Theorem \ref{general_single}}
	
	Start from
	\[
	\int_{\alpha}^{\beta}\int_{\alpha}^{\beta}\dfrac{\mathfrak{F}\left(zH\left(\vartheta_{1}\right)\right)-\mathfrak{F}\left(zH\left(\vartheta_{2}\right)\right)}{H\left(\vartheta_{1}\right)-H\left(\vartheta_{2}\right)}\,\mathrm{d}\vartheta_{1}\,\mathrm{d}\vartheta_{2}=\sum_{n=1}^{\infty}\mathfrak{R}_{n}z^{n}\int_{\alpha}^{\beta}\int_{\alpha}^{\beta}\frac{H^{n}\left(\vartheta_{1}\right)-H^{n}\left(\vartheta_{2}\right)}{H\left(\vartheta_{1}\right)-H\left(\vartheta_{2}\right)}\,\mathrm{d}\vartheta_{1}\,\mathrm{d}\vartheta_{2}
	.       \]
	Next, we define
	\[
	\int_{\alpha}^{\beta}H^{k}\left(\vartheta\right)\mathrm{d}\vartheta=\mathfrak{M}_{k},
	\]
	so that
	\[
	\int_{\alpha}^{\beta}\int_{\alpha}^{\beta}\frac{H^{n}\left(\vartheta_{1}\right)-H^{n}\left(\vartheta_{2}\right)}{H\left(\vartheta_{1}\right)-H\left(\vartheta_{2}\right)}\,\mathrm{d}\vartheta_{1}\,\mathrm{d}\vartheta_{2}=\sum_{k=0}^{n-1}\mathfrak{M}_{k}\mathfrak{M}_{n-k-1}=\mathfrak{B}_{n-1},
	\]
	is the coefficient of $z^{n-1}$ in the generating function
	\[
	G\left(z\right)=\left(\int_{\alpha}^{\beta}\frac{1}{1-zH\left(\vartheta\right)}\,\mathrm{d}\vartheta\right)^{2}=\sum_{k=1}^{\infty}\mathfrak{B}_{k}z^{k}=\left(\sum_{k= 0}^{\infty}\mathfrak{M}_{k}z^{k}\right)\left(\sum_{k= 0}^{\infty}\mathfrak{M}_{k}z^{k}\right).
	\]
	We deduce by the integral formula for the Hadamard product
	\begin{align*}
		\int_{\alpha}^{\beta}\int_{\alpha}^{\beta}\dfrac{\mathfrak{F}\left(zH\left(\vartheta_{1}\right)\right)-\mathfrak{F}\left(zH\left(\vartheta_{2}\right)\right)}{H\left(\vartheta_{1}\right)-H\left(\vartheta_{2}\right)}\,\mathrm{d}\vartheta_{1}\,\mathrm{d}\vartheta_{2} &=\sum_{n=1}^{\infty}\mathfrak{R}_{n}\mathfrak{B}_{n-1}\,z^{n}
		\\&  = \dfrac{1}{2\pi}\int_0^{2\pi}\mathfrak{F}\left(e^{\imath\vartheta}\sqrt{z}\right)e^{-\imath\vartheta}\sqrt{z}\,G\left(e^{\imath\vartheta}\sqrt{z}\right)\mathrm{d}\vartheta,
	\end{align*}
	which is the desired result. \QED
	\subsection{Proof of Theorem \ref{thm:Stieltjes}}
	
	Start from 
	\begin{align*}
		\int_{\mathbb{R}}\mathfrak{f}_{n+1}\left(x_{n+1}\right)\prod_{i=1}^{n}\hat{\mathfrak{f}}_{i}\left(x_{n+1}\right)\mathrm{d}x_{n+1} & =\int_{\mathbb{R}}\mathfrak{f}_{n+1}\left(x_{n+1}\right)\prod_{i=1}^{n}\int_{\mathbb{R}}\frac{\mathfrak{f}_{i}\left(x_{i}\right)}{1+x_{i}x_{n+1}}\,\mathrm{d}x_{i}\,\mathrm{d}x_{n+1}\\
		& =\int_{\mathbb{R}^{n}}\prod_{i=1}^{n}\mathfrak{f}_{i}\left(x_{i}\right)\mathrm{d}x_{i}\int_{\mathbb{R}}\frac{\mathfrak{f}_{n+1}\left(x_{n+1}\right)}{1+x_{i}x_{n+1}}\,\mathrm{d}x_{n+1}.
	\end{align*}
	Next, we apply the partial fraction decomposition to obtain
	\[
	\prod_{i=1}^{n}\frac{1}{1+x_{i}x_{n+1}}=\sum_{i=1}^{n}\frac{x_{i}^{n-1}}{\prod_{j\ne i}\left(x_{i}-x_{j}\right)},
	\]
	so that integrating on both sides yields
	\begin{align*}
		\int_{\mathbb{R}}\prod_{i=1}^{n}\frac{\mathfrak{f}_{n+1}\left(x_{n+1}\right)}{1+x_{i}x_{n+1}}\,\mathrm{d}x_{n+1} & =\sum_{i=1}^{n}\int_{\mathbb{R}}\frac{x_{i}^{n-1}}{\prod_{j\ne i}\left(x_{i}-x_{j}\right)}\mathfrak{f}_{n+1}\left(x_{n+1}\right)\mathrm{d}x_{n+1}\\
		& =\sum_{i=1}^{n}\frac{x_{i}^{n-1}}{\prod_{j\ne i}\left(x_{i}-x_{j}\right)}\hat{\mathfrak{f}}_{n+1}\left(x_{i}\right),
	\end{align*}
	and we deduce
	\[
	\int_{\mathbb{R}}\mathfrak{f}_{n+1}\left(x_{n+1}\right)\prod_{i=1}^{n}\hat{\mathfrak{f}}_{i}\left(x_{n+1}\right)\mathrm{d}x_{n+1}=\int_{\mathbb{R}^{n}}\prod_{i=1}^{n}\mathfrak{f}_{i}\left(x_{i}\right)\hat{\mathfrak{F}}_{n+1}\left(x_{1},x_2,\ldots,x_{n}\right)\mathrm{d}x_{i},
	\]
	as desired. \QED
	\subsection{Proof of Corollary \ref{cos_multivariant}}
	Using the partial fraction decomposition 
	\[
	\prod_{i=1}^{n}\dfrac{1}{1+\vartheta_{i}\vartheta_{n+1}}=\sum_{i=1}^{n}\dfrac{\vartheta_{i}^{n-1}}{\prod_{j=1,j\neq i}^{n}\left(\vartheta_{i}-\vartheta_{j}\right)}\dfrac{1}{1+\vartheta_{i}\vartheta_{n+1}},
	\]
	and the fact that 
	\[
	\int_{0}^{1}\dfrac{\vartheta_{i}^{n-1}}{1+\vartheta_{i}\vartheta_{n+1}}\,\mathrm{d}\vartheta_{n+1}=\vartheta_{i}^{n-2}\log\left(1+\vartheta_{i}\right),
	\]
	we have
	\begin{align*}
		&\int_{\left[0,1\right]^{n+1}}\prod_{i=1}^{n}\dfrac{1}{1+\vartheta_{i}\vartheta_{n+1}}\dfrac{\mathrm{d}\vartheta_{i}}{\sqrt{1-\vartheta_{i}^{2}}}\,\mathrm{d}\vartheta_{n+1}
		\\&=\int_{\left[0,1\right]^{n+1}}\sum_{i=1}^{n}\dfrac{\vartheta_{i}^{n-1}}{\prod_{j=1,j\neq i}^{n}\left(\vartheta_{i}-\vartheta_{j}\right)}\dfrac{1}{1+\vartheta_{i}\vartheta_{n+1}}\prod_{i=1}^{n}\dfrac{\mathrm{d}\vartheta_{i}}{\sqrt{1-\vartheta_{i}^{2}}}\,\mathrm{d}\vartheta_{n+1}
		\\&=\int_{\left[0,1\right]^{n+1}}\sum_{i=1}^{n}\dfrac{1}{\prod_{j=1,j\neq i}^{n}\left(\vartheta_{i}-\vartheta_{j}\right)}\dfrac{\vartheta_{i}^{n-1}}{1+\vartheta_{i}\vartheta_{n+1}}\prod_{i=1}^{n}\dfrac{\mathrm{d}\vartheta_{i}}{\sqrt{1-\vartheta_{i}^{2}}}\,\mathrm{d}\vartheta_{n+1}
		\\&=\int_{\left[0,1\right]^{n}}\sum_{i=1}^{n}\dfrac{\vartheta_{i}^{n-2}\log\left(1+\vartheta_{i}\right)}{\prod_{j=1,j\neq i}^{n}\left(\vartheta_{i}-\vartheta_{j}\right)}\prod_{i=1}^{n}\dfrac{\mathrm{d}\vartheta_{i}}{\sqrt{1-\vartheta_{i}^{2}}}
		\\&=\int_{\left[0,\pi/2\right]^{n}}\sum_{i=1}^{n}\frac{\cos^{n-2}\vartheta_{i}\log\cos\left(\vartheta_{i}/2\right)}{\prod_{j\ne i}\left(\cos\vartheta_{i}-\cos\vartheta_{j}\right)}\prod_{i=1}^{n}\mathrm{d}\vartheta_{i},
	\end{align*}
	where the last equality has simply followed by replacing $\vartheta_{i}\mapsto\cos \vartheta_{i}$.
	Next, we deduce 
	\[
	\int_{\left[0,\pi/2\right]^{n}}\sum_{i=1}^{n}\frac{\cos^{n-2}\vartheta_{i}\log\cos\left(\vartheta_{i}/2\right)}{\prod_{j\ne i}\left(\cos\vartheta_{i}-\cos\vartheta_{j}\right)}\prod_{i=1}^{n}\mathrm{d}\vartheta_{i}=\int_{0}^{1}\dfrac{\arccos^{n}x}{\left(1-x^{2}\right)^{\frac{n}{2}}}\,\mathrm{d}x=\int_{0}^{\pi/2}\dfrac{z^{n}}{\sin^{n-1}z}\,\mathrm{d}z,
	\]
	which is the desired result. \QED
	\subsection{Proof of Corollary \ref{tan_multivariant}}
	\label{proof of proposition 8.3}
	Start from the following identity
	\[
	\int_{0}^{\infty}\frac{1}{1+xz}\frac{\mathrm{d}x}{1+x^{2}}=\frac{1}{1+z^{2}}\left(\frac{\pi}{2}+z\log z\right).
	\]
	Now proceeding in a similar way as we did for proposition \ref{cos_multivariant}, we get
	\begin{align*}
		&\int_{\left[0,\infty\right)^{n}\times\left[0,1\right]}\prod_{i=1}^{n}\frac{1}{1+\vartheta_{i}\vartheta_{n+1}}\frac{\mathrm{d}\vartheta_{i}}{1+\vartheta_{i}^{2}}\,\mathrm{d}\vartheta_{n+1}\\&=\int_{\left[0,\infty\right)^{n}}\sum_{i=1}^{n}\dfrac{\vartheta_{i}^{n-2}\log\left(1+\vartheta_{i}\right)}{\prod_{j=1,j\neq i}^{n}\left(\vartheta_{i}-\vartheta_{j}\right)}\prod_{i=1}^{n}\dfrac{\mathrm{d}\vartheta_{i}}{1+\vartheta_{i}^{2}}
		\\& =\int_{\left[0,\pi/2\right)^{n}}\sum_{i=1}^{n}\frac{\tan^{n-2}\vartheta_{i}\log\tan\left(\vartheta_{i}/2\right)}{\prod_{j\ne i}\left(\tan\vartheta_{i}-\tan\vartheta_{j}\right)}\prod_{i=1}^{n}\mathrm{d}\vartheta_{i}\\&=\int_{0}^{1}\dfrac{1}{\left(1+z^{2}\right)^{n}}\left(\frac{\pi}{2}+z\log z\right)^{n}\mathrm{d}z,
	\end{align*}
	which is the desired result. \QED
	\section{Concluding Remarks and Further Research}
	The explicit computation of the main integral \eqref{double_integral}
	relies on the identity
	\[\int_{0}^{1}\dfrac{1}{1+zx}\frac{1}{\sqrt{1-x^2}}\, \mathrm{d}x = \dfrac{\arccos z}{\sqrt{1-z^2}},\]
	which expresses the (modified) Stieltjes transform of a centered Beta (or arcsine) distribution as a simple transformation of itself. This result is well-known in the resolvent method in the domain of random matrix theory. This suggests that integral \eqref{double_integral} could be related to probabilistic notions.
	Another path for future investigation is the link of this integral with MZVs as they appeared unexpectedly several times in the study.
	\section{Acknowledgements}
	The authors express their sincere gratitude to the anonymous reviewer for a very careful reading of the manuscript and many constructive comments and suggestions, which substantially helped us to improve the quality and expostion of the paper.
	
\end{document}